\documentclass{amsart}          

\usepackage{array}
\usepackage[margin=2.5cm]{geometry}
\usepackage[active]{srcltx}
\usepackage{amsfonts,amssymb,amsmath,latexsym,upref}
\usepackage[english]{babel}
\usepackage[latin1]{inputenc}  
\usepackage[all]{xy}\xyoption{arc}\xyoption{poly}
\usepackage{color}

\newtheorem{lemma}[equation]{Lemma}
\newtheorem{prop}[equation]{Proposition}
\newtheorem{thm}[equation]{Theorem}
\newtheorem{cor}[equation]{Corollary}

\newtheorem{defn}[equation]{Definition}

\theoremstyle{definition}
\newtheorem{exmp}[equation]{Example}

\newtheorem{rmk}[equation]{Remark}

\numberwithin{equation}{section}

\newcommand{\Z}{\mathbf{Z}}
\newcommand{\Q}{\mathbf{Q}}

\newcommand{\F}{\mathbf{F}}

\newcommand{\Aut}[2]{\operatorname{Aut}_{#1}(#2)}

\newcommand{\func}[3]{\mbox{$#1 \colon #2 \to #3$}}
\newcommand{\qt}[2]{#1 \backslash #2}

\newcommand{\GL}[3]{\operatorname{GL}^{#1}_{#2}(#3)}
\newcommand{\SL}[3]{\operatorname{SL}^{#1}_{#2}(#3)}
\newcommand{\PSL}[3]{\operatorname{PSL}^{#1}_{#2}(#3)}
\newcommand{\PGL}[3]{\operatorname{PGL}^{#1}_{#2}(#3)}

\newcommand{\Ob}[1]{\mathrm{Ob}(#1)}
\newcommand{\cat}[2]{\mathcal{#1}_{#2}^*} 
\newcommand{\catn}[3]{\mathcal{#1}_{#2}^{#3}}   
\newcommand{\catt}[2]{\widetilde{\mathcal{#1}}_{#2}^*} 
\newcommand{\catc}[2]{\mathcal{#1}_{#2}^c} 
\newcommand{\cattc}[2]{\widetilde{\mathcal{#1}}_{#2}^c} 

\newcommand{\m}{morphism}
\newcommand{\sfc}{$p$-selfcentralizing}

\newcommand{\syl}[1]{Sylow $#1$-subgroup}
\newcommand{\we}{weighting}
\newcommand{\Mb}{M\"obius}
\newcommand{\Euc}{Euler characteristic}

\title{Euler characteristics of $p$-subgroup categories}

\author{Martin Wedel Jacobsen}
\address{Institut for Matematiske Fag\\
  Universitetsparken 5\\
  DK--2100 K\o benhavn}
\email{tau.wedel@gmail.com}
\urladdr{}

\author{Jesper M.~M\o ller}
\address{Institut for Matematiske Fag\\
  Universitetsparken 5\\
  DK--2100 K\o benhavn}
\email{moller@math.ku.dk}
\urladdr{htpp://www.math.ku.dk/~moller}

\usepackage[bookmarks=true,bookmarksopen=false]{hyperref}
\hypersetup{
   pdftitle = {},
   pdfauthor = {Jesper Michael Møller},
   pdfpagemode = {UseOutlines},
   pdfstartview = {FitH},
   pdfborder = {0 0 0},
   backref = {true},
   colorlinks = {true},
   urlcolor = {blue},
   citecolor = {blue},
   linkcolor = {blue},
   pdftoolbar = {true},
}

 \keywords{\Euc, (co)\we, poset of
  $p$-subgroups, Frobenius category, \Mb\ function}

\begin{document}
\date{\today}
\maketitle
\tableofcontents

\begin{abstract}
  Let $G$ be a finite group and $p$ a prime number.  We compute the
  \Euc\ in the sense of Leinster for some categories of nonidentity
  $p$-subgroups of $G$. The $p$-subgroup categories considered include
  the poset $\catn SG*$, the transporter category $\catn TG*$, the
  linking category $\catn LG*$, the Frobenius, or fusion, category
  $\catn FG*$, and the orbit category $\catn OG*$ of all nonidentity
  $p$-subgroups of $G$.
\end{abstract}

\section{Introduction}
\label{sec:intro}

In this note we apply Tom Leinster's theory of \Euc s of (some) finite
categories \cite{leinster08} to $p$-subgroup categories associated to
finite groups.

For a finite group $G$ and a fixed prime number $p$, let
$\mathcal{S}_G$ denote the poset of consisting of all
$p$-subgroups of $G$ ordered by inclusion. Also, let $\mathcal{T}_G$
(the transporter category), $\mathcal{L}_G$ (the linking category
\cite{blo03}), $\mathcal{F}_G$ (the Frobenius category
\cite{puig09,blo03}), $\mathcal{O}_G$ (the orbit category), and
$\widetilde{\mathcal{F}}_G$ (the exterior quotient of the Frobenius
category \cite[1.3, 4.8]{puig09}) be the categories whose objects are
the $p$-subgroups of $G$ and whose \m\ sets are
\begin{alignat*}{3}
  &\mathcal{T}_G(H,K) = N_G(H,K) & &\qquad & 
  &\mathcal{L}_G(H,K) = \qt{O^pC_G(H)}{N_G(H,K)} \\
  &\mathcal{F}_G(H,K) = \qt{C_G(H)}{N_G(H,K)} & &\qquad & 
  &\mathcal{O}_G(H,K) =  N_G(H,K)/K  \\
  &\widetilde{\mathcal{F}}_G(H,K) = \qt{C_G(H)}{N_G(H,K)}/K
\end{alignat*}
for any two $p$-subgroups, $H$ and $K$, of $G$.  Here $N_G(H,K)=\{g
\in G \mid H^g \leq K \}$ denotes the transporter set.  Composition in
any of these categories is induced from group multiplication in $G$.
The \m s in $\mathcal{F}_G(H,K)$ are restrictions to $H$ of inner
auto\m s of $G$, \m s in $\mathcal{O}_G(H,K)$ are right $G$-maps $\qt
HG \to \qt KG$, and \m s in $\widetilde{\mathcal{F}}_G(H,K)$ are
$K$-conjugacy classes of restrictions to $H$ of inner auto\m s of $G$.
The endo\m\ groups in these categories of the $p$-subgroup $H$ of $G$
are $\catn SG*(H)=1$, $\catn TG{}(H) = N_G(H)$, $\catn LG{}(H) =
O^pC_G(H) \backslash N_G(H)$, $\catn FG{}(H) = C_G(H) \backslash
N_G(H)$, $\catn OG{}(H) = N_G(H)/H$, and $\widetilde{\mathcal{F}}_G(H)
= C_G(H) \backslash N_G(H) /H$.  The five categories $\mathcal{T}_G$,
$\mathcal{L}_G$, $\mathcal{F}_G$, $\mathcal{O}_G$, and
$\widetilde{\mathcal{F}}_G$ are related by functors
\begin{equation*}
  \xymatrix@R=10pt{
    {\mathcal{T}_G} \ar[r] \ar[d] &
    {\mathcal{L}_G} \ar[r] & 
    {\mathcal{F}_G} \ar[r] &
    {\widetilde{\mathcal{F}}_G} \\ 
    {\mathcal{O}_G} \ar[rrru] } 
\end{equation*}
If $\mathcal{C}$ is any of these categories
\begin{itemize}
\item $\mathcal{C}^*$ is the full subcategory of $\mathcal{C}$
  generated by all {\em nonidentity\/} $p$-subgroups
\item $\mathcal{C}^\mathrm{a}$ is the full subcategory of
  $\mathcal{C}^*$ generated by all {\em elementary abelian\/}
  $p$-subgroups
\item $\mathcal{C}^\mathrm{c}$ is the full subcategory of
  $\mathcal{C}^*$ generated by all {\em \sfc\/} $p$-subgroups
  (Definition~\ref{defn:sfc}) 
\end{itemize}

Here is a summary of our main results appearing in Table~\ref{tab:cat*},
Table~\ref{tab:wecat*}, and Theorem~\ref{thm:chiOG}.

\begin{thm}[\Euc s of $p$-subgroup categories]\label{thm:summary}
  The \Euc s are
  \begin{equation*}
    \chi(\mathcal{C}^*) =
    \sum_{[H]}\frac{-\mu(H)}{|\mathcal{C}^*(H)|}, \qquad 
    \mathcal{C} = \catn TG{}, \catn LG{}, \catn FG{}
  \end{equation*}
  where the sum runs over the set of conjugacy classes of nonidentity
  $p$-subgroups of $G$. Also, $\chi(\catn SG*) = |G|\chi(\catn TG*)$,
  $\chi(\widetilde{\mathcal{F}}_G^*) = \chi(\catn FG*)$, and
  \begin{equation*}
    \chi(\catn OG*) = 
        \chi(\catn TG*) + 
        \frac{p-1}{p} \sum_{\text{$[C]$ cyclic}} \frac{1}{|\catn
          OG*(C)|}
  \end{equation*}
  where the sum runs over the set of conjugacy classes of nonidentity
  cyclic $p$-subgroups of $G$.
  \end{thm}

  Here, $\mu(K)=\mu_K(1,K)$ is the \Mb\ function of finite groups
  \cite{HIO89} \cite[Chp 3.7]{stanley97}. (The formula for the \Euc\
  of $\cat SG$ is already known, of course.)
  Theorem~\ref{thm:summary} implies that
  \begin{equation*}
    \chi(\mathcal{C}^*) = \chi(\mathcal{C}^a), \qquad
    \mathcal{C} = \catn SG{}, \catn TG{}, \catn LG{}, \catn FG{} 
\end{equation*}
because these categories have co\we s supported (precisely) on the
elementary abelian subgroups of $G$. Dually, the \we s for $\catn
SG*$, $\catn TG*$, and $\catn OG*$ are supported on the nonidentity
$p$-radical subgroups of $G$ (Definition~\ref{defn:radical},
Corollary~\ref{cor:pradical}), and, if $G$ has a normal \syl p, $P$,
the \we\ for $\catn FG*$ is supported (precisely) on the subgroups of
$P$ of the form $C_P(x)$, $x \in G$ (Corollary~\ref{cor:normalsyl}).
This shows that these $p$-subgroup categories carry information,
retrieved by the \we\ or the co\we , about which objects are
elementary abelian, $p$-radical, or centralizers of group elements.
  The \Euc s of Theorem~\ref{thm:summary} are rational numbers.
  However, $|G|_{p'}\chi(\catn OG*)$ and $|G|_{p'}\chi(\catn FG*)$ are
  integers (Corollaries~\ref{cor:chiOGinZ} and \ref{cor:mchiGinZ}).

  As a spin-off of our investigations of (co)\we s we establish three
  combinatorial identities in \eqref{eq:combiSG}, \eqref{eq:combiFG},
  and \eqref{eq:combiOG}.

\begin{cor}\label{cor:combiids}
  For any finite group $G$ and any prime $p$,
  \begin{gather*}
    \sum_{H} \big(1-\chi(\catn S{\catn OG*(H)}*)+\mu(H) \big) = 0 \\
\sum_{H} \sum_{x \in C_G(H)} 
  \big( 1- \chi(\catn S{C_{\catn TG*(H)}(x)/H}*) 
  +\mu(H)\big) =0 \\
\sum_H \big( |H| - \chi(\catn S{\catn OG*(H)}*)|H|+\mu(H) \big) =
  \frac{p-1}{p}\sum_{C} |C|
  \end{gather*}
  where $H$ runs over the set of nonidentity $p$-subgroups of $G$ and $C$
  over the set of nonidentity cyclic $p$-subgroups of $G$.
\end{cor}

Theorems~\ref{thm:chiSprod} and \ref{thm:chiFprod} establish formulas
for \Euc s of posets and Frobenius categories of nonidentity subgroups
stating that
\begin{equation*}
   1-\chi(\cat S{\prod_{i=1}^n G_i}) = 
    \prod_{i=1}^n \big(1-\chi(\cat S{G_i})\big), \qquad
    1-\chi(\cat F{\prod_{i=1}^n G_i}) = 
    \prod_{i=1}^n \big(1-\chi(\cat F{G_i})\big) 
\end{equation*}
where $G_1,\ldots,G_n$ are finite groups.

For the sake of quick reference we list here the notation that we
are using throughout this paper:
\begin{itemize}
\item $p$ is a fixed prime number
\item $n_p$ is $p$-part of the integer $n$, the highest power of $p$
  dividing $n$, and $n_{p'} = n/n_p$ is the $p'$-part of $n$
\item $G$ is a finite group
\item $H \leq K$ means that $H$ is a subgroup of $K$ 
\item $\Phi(K)$ is the Frattini subgroup of $K$ \cite[Definition
3.14]{GLSII}
\item $\mathcal{C}$ is a finite category, $\mathcal{C}(a,b)$ is the
  set of \m s from object $a$ to object $b$, and $\mathcal{C}(a) =
  \mathcal{C}(a,a)$ is the monoid of endo\m s of $a$
\item $\Ob{\mathcal{C}}$ is the set of objects  of $\mathcal{C}$
\item $[\mathcal{C}]$ the set of iso\m\ classes of objects of
  $\mathcal{C}$ and $[a] \in [\mathcal{C}]$ the iso\m\ class of $a \in
  \Ob{\mathcal{C}}$
\end{itemize}



\section{\Euc s}
\label{sec:zetamu}

In this section we review the relevant parts of Tom Leinster's
concept of  \Euc\ of a finite category $\mathcal{C}$
\cite{leinster08}.

\subsection{The \Euc\ of a square matrix}
\label{sec:chisquare}

Let $S$ be a finite set and $\func{\zeta}{S \times S}{\Q}$ a rational
function on $S \times S$. Equivalently, $\zeta = \big( \zeta(a,b)
\big)_{a,b \in S}$ is a square matrix with rows and columns indexed by
the finite set $S$ and with rational entries $\zeta(a,b) \in \Q$, $a,b
\in S$.

\begin{defn}\cite[Definition 1.10]{leinster08}
  A \we\ for $\zeta$ is a column vector $(k^\bullet)$ and a
co\we\ for $\zeta$  is a row vector $(k_\bullet)$ solving the linear
equations 
\begin{equation*}
  \big(\zeta(a,b)\big)(k^b) =
  \begin{pmatrix}
    1 \\ \vdots \\ 1
  \end{pmatrix}, \qquad
  (k_a)\big(\zeta(a,b)\big) =
  \begin{pmatrix}
    1 & \cdots & 1
  \end{pmatrix}
\end{equation*}
\end{defn}

If $\zeta$ admits both a \we\ $k^\bullet$ and a co\we\ 
$k_\bullet$ then the sum of the values of the \we\
\begin{equation}
  \label{eq:chi}
   \sum_{b} k^b =
    \sum_{b} 
    \big(\sum_{a} k_a\zeta(a,b))\big) k^b =
    \sum_{a} k_a 
    \big(\sum_{b} \zeta(a,b)k^b \big) =
     \sum_{a} k_a
\end{equation}
equals the sum of the values of the co\we\ (and then this sum is
independent of the choice of \we\ or co\we ).

\begin{defn}\cite[Definition 2.2]{leinster08}\label{defn:chi}
  The square matrix $\zeta$ {\em has \Euc\/} if it admits both
  a \we\  and a co\we\  and its \Euc\ is then
  the sum
  \begin{equation*}
    \chi(\zeta) =  \sum_{b} k^b = \sum_{a} k_a
  \end{equation*}
  of all the values of any \we\ $k^\bullet$ or any co\we\ $k_\bullet$.
\end{defn}

In case the square matrix $\zeta$ is invertible, if we let $\mu=
\big( \mu(a,b) \big)_{a,b \in S}$ denote the inverse of $\zeta$, the
{\em \Mb\ inversion formula\/}
\begin{equation}\label{eq:zetamu}
  \forall a,c \in S \colon
   \sum_b \zeta(a,b)\mu(b,c) = \delta(a,c) = \sum_b \mu(a,b)\zeta(b,c) 
\end{equation}
simply expresses that $\zeta$ and $\mu$ are inverse matrices. When
$\zeta$ is invertible the vectors
\begin{equation*}
  \big( k^a \big) = \big( \mu(a,b) \big)
  \begin{pmatrix}
    1 \\ \vdots \\ 1
  \end{pmatrix} = \big( \sum_{b \in S} \mu(a,b) \big),
  \qquad 
  \big( k_b \big) =
  \begin{pmatrix}
    1 & \hdots & 1
  \end{pmatrix}
  \big( \mu(a,b) \big) 
  = \big( \sum_{a \in S} \mu(a,b) \big)
\end{equation*}
are the {\em unique\/} \we\ and co\we\ for $\zeta$ and the \Euc\ of
$\zeta$
\begin{equation*}
  \chi(\zeta) = \sum_{a,b \in S} \mu(a,b)
\end{equation*}
is the sum of all the entries in the inverse matrix.  

\subsection{The \Euc\ of a finite category}
\label{sec:chifinitecat}
Define the $\zeta$-matrix for the finite category $\mathcal{C}$ to be
the square matrix
\begin{equation*}
  \zeta(\mathcal{C}) = 
  \big( |\mathcal{C}(a,b)| \big)_{a,b \in \Ob{\mathcal{C}}}
\end{equation*}
that counts the number of \m s between pairs of objects of
$\mathcal{C}$.  We say that the category $\mathcal{C}$ admits a \we ,
admits a co\we , or has \Euc\ if its $\zeta$-matrix
$\zeta(\mathcal{C})$ does.  This means that a \we\ for $\mathcal{C}$
is a rational function \func{k^\bullet}{\Ob{\mathcal{C}}}{\Q} and a
co\we\ for $\mathcal{C}$ is a rational function
\func{k_\bullet}{\Ob{\mathcal{C}}}{\Q} such that
    \begin{equation}
      \label{eq:weight} 
      \forall a \in \Ob{\mathcal{C}} \colon 
      \sum_{b \in \Ob{\mathcal{C}}} \zeta(a,b)k^b = 1, \qquad
      \forall b \in \Ob{\mathcal{C}} \colon 
      \sum_{a \in \Ob{\mathcal{C}}} k_a\zeta(a,b) = 1
    \end{equation}
and the \Euc\ of $\mathcal{C}$ is
\begin{equation*}
  \chi(\mathcal{C}) = \sum_{b \in \Ob{\mathcal{C}}}k^b
  = \sum_{a \in \Ob{\mathcal{C}}}k_a
\end{equation*}
provided that $\mathcal{C}$ admits both a \we\ and a co\we .  We say
that $\mathcal{C}$ has \Mb\ inversion if $\zeta(\mathcal{C})$ is
invertible and then
\begin{equation*}
  k^a = \sum_{b \in \Ob{\mathcal{C}}} \mu(a,b), \qquad
  k_b = \sum_{a \in \Ob{\mathcal{C}}} \mu(a,b), \qquad
  \chi(\mathcal{C}) = \sum_{a,b \in \Ob{\mathcal{C}}} \mu(a,b)
\end{equation*}
is the unique \we , the unique co\we , and the \Euc\ of $\mathcal{C}$
where $\mu=\zeta(\mathcal{C})^{-1}$ denotes the inverse of the
$\zeta$-matrix.


\begin{exmp}\label{exmp:finalelement}\cite[Examples 1.1.c]{leinster08}
  Suppose that $\mathcal{C}$ has \Euc . If $\mathcal{C}$ has an
  initial element $0$ then the Kronecker function $k_\bullet =
  \delta(0,\bullet)$ is a
  co\we\  concentrated at the initial element because
  \begin{equation*}
    \sum_a \delta(0,a)\zeta(a,b) = \zeta(0,b) = 1
  \end{equation*}
  and the \Euc\ of $\mathcal{C}$ is $\chi(\mathcal{C}) = \sum_a
  \delta(0,a) = \delta(0,0)= 1$.  Dually, if $\mathcal{C}$ has a
  terminal element $1$, then $k^\bullet = \delta(\bullet,1)$ is a
  weighting concentrated at the terminal element, and
  $\chi(\mathcal{C})=1$.
\end{exmp}

As usual, $\delta$ stands for Kronecker's $\delta$-function
\begin{equation*}
  \delta(a,b) =
  \begin{cases}
    1 & a = b \\ 0 & a \neq b
  \end{cases},  \qquad a,b \in \Ob{\mathcal{C}},
\end{equation*}

\begin{lemma}\cite[Proposition 2.4]{leinster08}\label{lemma:adjointcats}
  Let $\mathcal{C}$ and $\mathcal{D}$ be finite categories.
  \begin{enumerate}
  \item $\mathcal{C}$ has \Euc\ if and only if its opposite category
    $\mathcal{C}^{\mathrm{op}}$ has and then $\chi(\mathcal{C}) =
    \chi(\mathcal{C}^{\mathrm{op}})$.
  \item If both $\mathcal{C}$ and $\mathcal{D}$ have \Euc s and there
    is an adjunction $\xymatrix@1{{\mathcal{C}} \ar@<.5ex>[r] &
      {\mathcal{D}} \ar@<.5ex>[l]}$ then
    $\chi(\mathcal{C})=\chi(\mathcal{D})$.
  \item If $\mathcal{C}$ and $\mathcal{D}$ are equivalent then
    $\mathcal{C}$ has \Euc\ if and only if $\mathcal{D}$ has \Euc\ and
    then $\chi(\mathcal{C})=\chi(\mathcal{D})$.
  \end{enumerate}
\end{lemma}

\begin{lemma}\label{lemma:supportk}
  Let $\mathcal{C}$ be a full subcategory of $\mathcal{D}$ and suppose
  that both categories have \Euc s.
  \begin{enumerate}
  \item If $\Ob{\mathcal{C}}$ contains the support of some \we\
    $k^\bullet$ on $\mathcal{D}$, then the restriction $k^\bullet
    \vert \Ob{\mathcal{C}}$ is a \we\ for $\mathcal{C}$ and
    $\chi(\mathcal{C}) = \chi(\mathcal{D})$.\label{lemma:supportk1}
  \item If $\Ob{\mathcal{C}}$ contains the support of some co\we\
    $k_\bullet$ on $\mathcal{D}$, then the restriction $k_\bullet
    \vert \Ob{\mathcal{C}}$ is a co\we\ for $\mathcal{C}$ and
    $\chi(\mathcal{C}) = \chi(\mathcal{D})$.\label{lemma:supportk2}
  \end{enumerate}
\end{lemma}
\begin{proof}
  \eqref{lemma:supportk2} The assumption is that $\forall a \in
  \Ob{\mathcal{D}} \colon k_a \neq 0 \Longrightarrow a \in
  \Ob{\mathcal{C}}$. For any $b \in \Ob{\mathcal{C}}$
  \begin{equation*}
    \sum_{a \in \Ob{\mathcal{C}}} k_a \zeta(a,b) = 
     \sum_{a \in \Ob{\mathcal{D}}} k_a \zeta(a,b) = 1
  \end{equation*}
  so that the restriction to $\Ob{\mathcal{C}}$ of $k_\bullet$ is
  indeed a co\we\  for $\mathcal{C}$. The \Euc\ of $\mathcal{C}$
  is $\chi(\mathcal{C}) = \sum_{a \in \Ob{\mathcal{C}}}k_a = \sum_{a
    \in \Ob{\mathcal{D}}}k_a = \chi(\mathcal{D})$.
\end{proof}


\subsection{The \Euc\ of a finite poset}
\label{sec:chiposet}

In particular, any finite {\em poset}, $\mathcal{S}$, has \Mb\
inversion \cite{stanley97}. The value of the \Mb\ function
\begin{equation*}
  \mu(a,b) = \chi((a,b))-1, \qquad a<b,
\end{equation*}
depends only on the open interval $(a,b)$ from $a$ to $b$ and not on
the whole poset \cite[Proposition 3.8.5]{stanley97} \cite[Corollary
1.5]{leinster08}.

\begin{exmp}\label{exmp:S*}
  Let $\mathcal{S}$ be a finite poset with a least element, $0$.  Then
  $\chi(\mathcal{S})=1$ by Example~\ref{exmp:finalelement}.  For any
  element $b$ of $\mathcal{S}$,
\begin{equation*}
  \sum_{a \colon 0 \leq a \leq b}\mu(a,b) 
  =  \sum_{a \colon 0 \leq a \leq b}|\mathcal{S}(0,a)| \mu(a,b)
  =  \delta(0,b)
\end{equation*}
and for any element $b>0$ of $\mathcal{S}$,
\begin{equation}
  \label{eq:sumSneq0}
  \sum_{0 \lneqq a \leq b}\mu(a,b) 
  = -\mu(0,b) +   \sum_{0 \leq a \leq b}\mu(a,b) = -\mu(0,b)
\end{equation}
The \Euc\ of the subposet $\mathcal{S}^*$ of all elements $\neq 0$ is
\begin{equation*}
  \chi(\mathcal{S}^*) = \sum_{a,b>0} \mu(a,b) = \sum_{b>0} \sum_{a>0}
  \mu(a,b) = \sum_{b>0} -\mu(0,b)
\end{equation*}
Alternatively,
\begin{equation}
  \label{eq:altmuposet}
  1-\chi(\mathcal{S}^*) = \mu(0,0) + \sum_{b>0}\mu(0,b) = 
  \sum_{b \in \Ob{\mathcal{S}}}\mu(0,b)
\end{equation}
with summation over {\em all\/} elements of the poset $\mathcal{S}$. 
\end{exmp}

\subsection{The \Euc\ of $[\mathcal C]$}
\label{sec:pi0C}

Let $[\func{\zeta(\mathcal{C})]}{ [\mathcal{C}]\times
  [\mathcal{C}]}{\Q}$ be the function induced by the $\zeta$-function
$\func{\zeta(\mathcal{C})}{\Ob{\mathcal{C}}\times\Ob{\mathcal{C}}}{\Q}$
for $\mathcal{C}$.  We say that the set $[\mathcal{C}]$ of iso\m\
classes of $\mathcal{C}$-objects admits a \we , a co\we , or has \Euc\
if its $\zeta$-matrix $[\zeta(\mathcal{C})]$ does.  This means that a
\we\ for $ [\mathcal{C}]$ is a rational function \func{k^\bullet}{
  [\mathcal{C}]}{\Q} and a co\we\ for $ [\mathcal{C}]$ is a rational
function \func{k_\bullet}{ [\mathcal{C}]}{\Q} such that
    \begin{equation}
      \label{eq:pi0weight}
      \forall [a] \in  [\mathcal{C}] \colon 
      \sum_{[b] \in  [\mathcal{C}]} [\zeta(\mathcal{C})]([a],[b])k^{[b]} = 1,
      \qquad 
      \forall [b] \in  [\mathcal{C}] \colon 
      \sum_{[a] \in  [\mathcal{C}]} k_{[a]}[\zeta(\mathcal{C})]([a],[b]) = 1
    \end{equation}
    and the \Euc\ of $ [\mathcal{C}]$ is
\begin{equation*}
  \chi( [\mathcal{C}]) = \sum_{[b] \in  [\mathcal{C}]} k^{[b]}
  = \sum_{[a] \in  [\mathcal{C}]} k_{[a]}
\end{equation*}
provided that $ [\mathcal{C}]$ admits both a \we\ and a co\we .
Clearly, if $\mathcal{C}$ has a \we\ $k^\bullet$ and a co\we\
$k_\bullet$, then $ [\mathcal{C}]$ has \we\ $k^{[a]} = \sum_{c \in
  [a]}k^c$ and co\we\  $k_{[b]} = \sum_{c \in [b]}k_c$ and
$\chi(\mathcal{C}) = \chi( [\mathcal{C}])$.

We say that $ [\mathcal{C}]$ has \Mb\ inversion if its $\zeta$-matrix
$[\zeta(\mathcal{C})]$ is invertible and then
\begin{equation*}
  k^{[a]} = \sum_{[b] \in  [\mathcal{C}]} [\mu]([a],[b]), \qquad
  k_{[b]} = \sum_{[a] \in  [\mathcal{C}]} [\mu]([a],[b]), \qquad
  \chi( [\mathcal{C}]) = \sum_{[a],[b] \in  [\mathcal{C}]}
  [\mu]([a],[b]) 
\end{equation*}
is the unique \we , the unique co\we , and the \Euc\ of $
[\mathcal{C}]$ where $\big( [\mu]([a],[b] \big)_{[a],[b] \in
  [\mathcal{C}]}$ denotes the inverse of $[\zeta(\mathcal{C})]$.

\begin{thm}\label{thm:wefrompi0C}
    Suppose that $ [\mathcal{C}]$ has \Mb\ inversion. 
Then the functions
\begin{equation*}
  k^a = |[a]|^{-1} k^{[a]}, \qquad
  k_b = |[b]|^{-1} k_{[b]}
\end{equation*}
are a \we\ and a co\we\ for $\mathcal{C}$, and
$\chi( [\mathcal{C}]) = \chi(\mathcal{C})$.
\end{thm}
\begin{proof}
  Let $a$ be any object of $\mathcal{C}$. Since the function
  $k^\bullet$ is constant on the iso\m\ class $[a]$ we find that
  \begin{equation*}
    \sum_b \zeta(a,b)k^b =
    \sum_b [\zeta]([a],[b])k^b =
    \sum_{[b]} [\zeta]([a],[b])k^b |[b]| =
    \sum_{[b]} [\zeta]([a],[b]) k^{[b]} =
    1
  \end{equation*}
  according to \eqref{eq:pi0weight} and this shows that $k^\bullet$ is
  a \we\ on $\mathcal{C}$ according to \eqref{eq:weight}. A symmetric
  argument shows that $k_\bullet$ is a co\we . Thus $\mathcal{C}$ has
  \Euc\ $\chi(\mathcal{C}) = \sum_a k^a = \sum_{[a]} |[a]|k^a =
  \sum_{[a]}k^{[a]} = \chi( [\mathcal{C}])$.
\end{proof}

Observe for instance that the transporter category $\catn TG*$ (in
general) does not have \Mb\ inversion as its $\zeta$-matrix
$\zeta(\catn TG*)$ has two identical rows as soon as there are two
nonidentical nonidentity subgroups of $G$ that are conjugate in $G$.
However, $[\catn TG*]$, the set of conjugacy classes of subgroups of
$G$, always has \Mb\ inversion: Extend the partial subconjugation
ordering, $[H] \leq_G [K] \iff \catn TG*(H,K) \neq \emptyset$, to a
total ordering of $[\catn TG*]$. If $[H]>[K]$ in this total order,
then $\catn TG*(H,K) = \emptyset$. This means that the $\zeta$-matrix
$(|\catn TG*(H,K)|)_{[H],[K] \in [\catn TG*]}$ is upper-triangular in
this total ordering and, as the diagonal entries are nonzero, it is
invertible.  We shall determine the \Mb\ function $[\mu]$ for $[\catn
TG*]$ in Proposition~\ref{prop:nupi0TG*}.

\subsection{The \Euc\ of a  homotopy orbit category}
\label{sec:ShG}

Let $ {\mathcal{S}}$ be a finite category with a $G$-action.  (This
means that there is a functor from $G$ to the category of finite
categories taking the single object of $G$ to $\mathcal{S}$.)  The
{\em homotopy orbit category\/}, $ {\mathcal{S}}_{hG}$, is the
Grothendieck construction on the $G$-action on $ {\mathcal{S}}$: The
category with the same set of objects as $ \mathcal{S}$ and with \m\
sets
\begin{equation}\label{eq:ShG}
  {\mathcal{S}}_{hG}(a,b) = \coprod_{g \in G}
  {S}(ag,b)
  = \coprod_{g \in G} {S}(a,bg^{-1}) \qquad\text{and}\qquad 
  |{\mathcal{S}}_{hG}(a,b)| 
  = \sum_{g \in G} |{S}(a,bg^{-1})| 
\end{equation}


\begin{thm}\label{thm:ShG}
  Let $\mathcal{F}$ be a finite category with the same objects as $ S$
  such that $d_a |\mathcal{F}(a,b)| t^b = |\mathcal{S}_{hG}(a,b)|$ for
  all $a,b \in \Ob{ S}$.
  \begin{enumerate}
  \item If $m^\bullet \colon \Ob{\mathcal{S}} \to \Q$ is rational
    function so that $\sum_b |\mathcal{S}(a,b)|m^b = d_a$ for all $a
    \in \Ob{\mathcal{S}}$ and $d_\bullet$ is $G$-invariant, then
    $|G|^{-1}t^\bullet m^\bullet$ is a \we\ for $\mathcal{F}$
    \label{prop:cowetildeF1}
  \item If $m_\bullet \colon \Ob{\mathcal{S}} \to \Q$ is rational
    function so that $\sum_a m_a |\mathcal{S}(a,b)| = t^b$ for all $b
    \in \Ob{\mathcal{S}}$ and $t^\bullet$ is $G$-invariant, then
    $|G|^{-1} m_\bullet d_\bullet$ is a co\we\ for $\mathcal{F}$
    \label{prop:cowetildeF2}
  \item Suppose that $\mathcal{S}$ has \Mb\ inversion and $\mu$ is the
    \Mb\ function. If $d_\bullet$ is $G$-invariant then $k^a =
    |G|^{-1} \sum_{b} t^a \mu(a,b) d_b$ is a \we\ for $\mathcal{F}$,
    and if $t^\bullet$ is $G$-invariant then $k_b = |G|^{-1} \sum_{a}
    t^a \mu(a,b) d_b$ a co\we\ for $\mathcal{F}$.
    \label{prop:cowetildeF3} 
  \end{enumerate}
\end{thm}
\begin{proof}
  \noindent \eqref{prop:cowetildeF1} The proofs of
  \eqref{prop:cowetildeF1} and \eqref{prop:cowetildeF2} are dual to
  each other.

  \noindent \eqref{prop:cowetildeF2}
  For every $b \in \Ob{\mathcal{S}}$,
  \begin{equation*}
    \sum_a m_a d_a |\mathcal{F}(a,b)|
    = \sum_a m_a |\mathcal{S}_{hG}(a,b)|(t^b)^{-1}
    =  \sum_{g \in G} \sum_a m_a |\mathcal{S}(a,bg^{-1})|(t^b)^{-1}
    =  \sum_{g \in G} t^{bg^{-1}}(t^b)^{-1} 
    = |G| 
  \end{equation*}
  as $t^\bullet$ is $G$-invariant so that $t^{bg^{-1}} = t^b$ for all
  $g \in G$.

  \noindent \eqref{prop:cowetildeF3}  If $m^a = \sum_b \mu(a,b)d_b$
  then $\sum_b |\mathcal{S}(a,b)|m^b = d_a$ by the \Mb\ inversion
  formula \eqref{eq:zetamu}. By \eqref{prop:cowetildeF1}, $k^a$ is a
  \we\ for $\mathcal{F}$ if $d_\bullet$ is $G$-invariant.
  Dually,
  If $m_b = \sum_a t^a\mu(a,b)$
  then $\sum_a m_a |\mathcal{S}(a,b)| = t^b$ by the \Mb\ inversion
  formula \eqref{eq:zetamu}. By \eqref{prop:cowetildeF2}, $k_b$ is a
  co\we\ for $\mathcal{F}$ if $t^\bullet$ is $G$-invariant.
\end{proof}



\section{The \Mb\ function of a finite group}
\label{sec:chiSG}


The \Mb\ function for the finite group $G$ is the \Mb\ function $\mu$
for the poset $\overline{\mathcal{S}}_G$ of {\em all\/} subgroups of $G$.
Note that $\mu$ restricts to \Mb\ functions for the convex subposets
$\catn SG{}$ and $\catn SG*$ of $p$-subgroups.  For any subgroup $K
\leq G$, $\mu(1,K)$ only depends on $K$ and not on the whole group $G$
and it is customary to write $\mu(K)$ for $\mu(1,K)$ \cite{HIO89}.

\begin{lemma}\cite{hall36}\cite[Corollary 3.5]{HIO89}
  \label{lemma:muH}
  Let $H$ and $K$ be $p$-subgroups of $G$. Then
  \begin{equation*}
    \mu(H,K) =
    \begin{cases}
       (-1)^np^{\binom{n}{2}} & \Phi(K) \leq H \leq K, \quad p^n=|K
       \colon H| \\
       0 & \text{otherwise}
    \end{cases}
  \end{equation*}
  In particular, 
  $\mu(K)=\mu(1,K)=0$ unless $K$ is
  elementary abelian where
  \begin{equation*}
    \mu(K) = (-1)^np^{\binom{n}{2}}, \quad p^n = |K|
  \end{equation*}
\end{lemma}
\begin{proof} 
  If $\mu(H,K) \neq 0$ then $H \lhd K$ with $H \backslash K$
  elementary abelian and $\mu(H,K) = \mu(H \backslash K)$
  \cite[Proposition 2.4]{kratzer_thevenaz84} \cite[Lemme
  4.1]{kratzer_thevenaz85}. Burnside's basis theorem
  \cite[5.3.2]{robinson:groups} \cite[Lemma 3.15]{GLSII}, $\Phi(K) =
  [K,K]K^p$, shows that $H \lhd K$ with $H \backslash K$ elementary
  abelian if and only if $\Phi(K) \leq H$.
\end{proof}

\begin{table}[t]
  \centering
  \begin{tabular}{c||l|l|l}
    $ \mathcal{C}$ &  $k^H$  &  $k_K$  &  $\chi(\mathcal{C})$  \\ \hline
    $\displaystyle \catn SG*$ &  $\displaystyle \sum_K \mu(H,K)$
    &  $\displaystyle-\mu(K)$  &  
    $\displaystyle     |G| \chi(\catn TG*)$\\ 
    $\displaystyle \catn TG*$   &  
    $\displaystyle|G|^{-1}\sum_K \mu(H,K)$  & 
    $\displaystyle -|G|^{-1}\mu(K)$  &  
    $\displaystyle \sum_{[K]} \frac{-\mu(K)}{|\catn TG*(K)|}$  \\
    $\catn LG*$  &  $\displaystyle|G|^{-1}\sum_K \mu(H,K) |O^pC_G(K)|$  &
    $\displaystyle -|G|^{-1}\mu(K)|O^pC_G(K)|$  &
    $\displaystyle  \sum_{[K]} \frac{-\mu(K)}{|\catn LG*(K)|}$  \\
    $\displaystyle \catn FG*$  &  
    $\displaystyle|G|^{-1}\sum_K \mu(H,K) |C_G(K)|$  &
    $\displaystyle-|G|^{-1}\mu(K)|C_G(K)|$  &
    $\displaystyle \sum_{[K]} \frac{-\mu(K)}{|\catn FG*(K)|}$  \\
    $\displaystyle \catn OG*$  &  
    $\displaystyle|G|^{-1}|H|\sum_K\mu(H,K)$  &
    $\displaystyle |G|^{-1}\sum_H |H| \mu(H,K)$  &  
    $\displaystyle |G|^{-1}\sum_{H,K} |H| \mu(H,K)$ 
  \end{tabular}
  \caption{Categories of nonidentity
    $p$-subgroups}
  \label{tab:cat*}
\end{table}


\begin{thm}\label{thm:STLFO}
  Weightings $k^\bullet$, co\we s $k_\bullet$, and \Euc s for the
  $p$-subgroup categories $\mathcal{S}_G^*$, $\mathcal{T}_G^*$,
  $\mathcal{L}_G^*$, $\mathcal{F}_G^*$, and $\mathcal{O}_G^*$ are as
  in Table~\ref{tab:cat*}.
\end{thm}
\begin{proof}
  This follows almost immediately from Theorem~\ref{thm:ShG} because
  the transporter category $\cat TG = (\cat SG)_{hG}$ is the homotopy
  orbit category for the conjugation action of $G$ on the poset $\cat
  SG$ and
  \begin{equation*}
    |\cat TG(H,K)| = |O^pC_G(H)| |\cat LG(H,K)| = |C_G(H)| |\cat
  FG(H,K)| = |\cat OG(H,K)| |K|
  \end{equation*}  
  Since the poset $\catn SG*$ has \Mb\ inversion,
  $k^H=|G|^{-1}\sum_K\mu(H,K)$ is a \we ,$k_K = |G|^{-1}\sum_H\mu(H,K)
  = -|G|^{-1}\mu(K)$ (Example~\ref{exmp:S*}) a co\we\ for $\catn TG*$
  by Theorem~\ref{thm:ShG}.\eqref{prop:cowetildeF3}.  In case of
  $\catn FG*$, Theorem~\ref{thm:ShG}.\eqref{prop:cowetildeF3} provides
  the \we\ and the co\we\
  \begin{equation*}
    k^H = |G|^{-1} \sum_{H \leq K} \mu(H,K) |C_G(K)|, \qquad
    k_K = |G|^{-1} \sum_{H \leq K} \mu(H,K) |C_G(K)| =  -|G|^{-1} \mu(K)
    |C_G(K)|  
  \end{equation*}
  where the expression for the co\we\ simplifies when using $\sum_{H
    \in \Ob{\mathcal{S}_K^*}} \mu(H,K) = -\mu(K)$ from identity
  \eqref{eq:sumSneq0}. The \Euc\ of $\catn FG*$, calculated as the sum
  of the values of the co\we , is
  \begin{equation*}
    \chi(\catn FG*) = \sum_K k_K 
                    = -|G|^{-1} \sum_K  \mu(K) |C_G(K)|
                    = -|G|^{-1} \sum_{[K]} \mu(K) |C_G(K)| |G \colon
                    N_G(K)| 
                    = \sum_{[K]} \frac{-\mu(K)}{|\catn FG*(K)|} 
  \end{equation*}
  because the co\we\ $k_K$ is constant over the conjugacy class $[K]$
  of $K$ and $|[K]| = |G \colon N_G(K)|$.
\end{proof}

The quotient category $\catt FG$ is missing from Table~\ref{tab:cat*}
because Theorem~\ref{thm:ShG} does not directly apply. We shall later
see that $\catt FG$ and $\catn FG*$ have identical \Euc s
(Corollary~\ref{cor:chiFGeqchiFtildeG}).

Lemma~\ref{lemma:supportk} implies that $\chi(\mathcal{C}^*) =
\chi(\mathcal{C}^a)$ for $\mathcal{C}= \mathcal{S}, \mathcal{T},
\mathcal{L}, \mathcal{F}$ because the co\we s for these categories are
concentrated on the elementary abelian $p$-subgroups of $G$
(Lemma~\ref{lemma:muH}).  Quillen shows in \cite[Proposition
2.1]{quillen78} the much stronger result that the posets
$\mathcal{S}_G^*$ and $\mathcal{S}_G^{\mathrm{a}}$ are homotopy
equivalent.


If $P$ is a nonidentity $p$-group we immediately have that
\begin{equation}
  \label{eq:chiSP}
  \chi(\cat SP) = 1, \quad 
  \chi(\cat TP) = |P|^{-1}, \quad
  \chi(\cat LP) = |P|^{-1}, \quad
  \chi(\cat FP) = 1, \quad
  \chi(\cat OP) = 1, \quad 
  \chi(\catt FP) = 1
\end{equation}
because $P$ is terminal in $\cat SP$ and $\cat OP$, $\catn TP* = \catn
LP*$ and $\chi(\catn TP*) = |P|^{-1} \chi(\catn SG*) = |P|^{-1}$ by
Theorem~\ref{thm:ShG}, Proposition~\ref{prop:chiFG=1} applies to $\cat
FP$, and Corollary~\ref{cor:chiFGeqchiFtildeG} to $\catt FP$.  More
generally, if $G$ has a normal $p$-complement, then $\chi(\cat FG)=1$
because $\cat FG = \cat FP$ according to the Frobenius normal
$p$-complement theorem \cite[Proposition
16.10]{GLSIII}\cite[10.3.2]{robinson:groups}.  For some more examples,
let $D_{pn}$ be the dihedral group of order $2pn$, $n \geq 1$, $A_p$
the alternating group of order $p>2$, and $\SL{}n{\F_q}$ the special
linear group where $q$ is a power of $p$ and $n \geq 2$.  Then
  \begin{equation*}
    \chi(\cat S{D_{pn}}) = 1, \qquad
    \chi(\cat S{A_p}) = (p-2)!, \qquad
    \chi(\cat S{\SL{}n{\F_q}}) = 1+ (-1)^n q^{\binom n2}
  \end{equation*}
  See \cite[Example 2.7]{quillen78} for the \Euc\ of $\cat S{A_p}$.
  Let $V_n(q)$ be an $n$-dimensional vector space over $\F_q$ and
  $L_n(q)$ the poset of $\F_q$-subspaces of $V_n(q)$. As $\mathcal{S}
  _{\SL{}n{\F_q}}^{\mathrm{a}})$ and the open interval $(0,V_n(q))$,
  the building for $\SL{}n{\F_q}$ \cite[Example
  6.5]{abramenko_brown2008}, are homotopy equivalent posets
  \cite[Theorem 3.1]{quillen78},
  \begin{equation*}
    \chi(\mathcal{S}_{\SL{}n{\F_q}}^{\mathrm{a}}) 
    = \chi((0,V_n(q))) 
    = 1+\mu_{L_n(q)}(0,V_n(q)) 
    = 1+ (-1)^n q^{\binom n2}
  \end{equation*}
  by the computation of the \Mb\ function $\mu_{L_n(q)}$ in $L_n(q)$
  \cite[Example 3.10.2]{stanley97} \cite[Proposition
  3.6]{kratzer_thevenaz85}.  In this example we may replace
  $\SL{}n{\F_q}$ by any of the groups $\GL{}n{\F_q}$, $\PSL{}n{\F_q}$,
  or $\PGL{}n{\F_q}$ since they all have identical $p$-subgroup
  posets.  The computer-generated Table~\ref{tab:FAn} displays \Euc s
  of poset categories at $p=2$ of small alternating groups.

 

  The \Euc s of the subgroup categories
  generated by {\em all\/} $p$-subgroups of $G$ (including the
  identity subgroup) are
  \begin{multline*}
    \chi(\mathcal{S}_G)=1, \quad
    \chi(\mathcal{T}_G)=|G|^{-1}, \quad
    \chi(\mathcal{L}_G)=|G \colon O^pG|^{-1}, \quad
    \chi(\mathcal{F}_G)=1, \\
    \chi(\mathcal{O}_G)=
    |G|^{-1} + \frac{p-1}{p} \sum_{[C]}\frac{1}{|\cat OG(C)|}, \quad
    \chi(\widetilde{\mathcal{F}}_G)=1
  \end{multline*}
  Observe that $\mathcal{S}_G$, $\mathcal{F}_G$, and
  $\widetilde{\mathcal{F}}_G$ have initial objects and that
  $\mathcal{T}_G$ deformation retracts onto $\mathcal{T}_G(1) = G$ and
  $\mathcal{L}_G$ deformation retracts onto $\mathcal{L}_G(1) =
  \qt{O^pG}G$. See Remark~\ref{rmk:chiOG} for $\chi(\mathcal{O}_G)$.



\subsection{\Euc\ of the exterior quotient of the Frobenius category}
\label{sec::extquotFG}

The equation
\begin{equation}\label{eq:tildeFG1}
   |C_G(H)| |\catt FG(H,K)| |K| =  \sum_{n \in \cat TG(H,K)}
   |C_K(H^n)|  
\end{equation}
follows from Burnside's counting lemma (Lemma~\ref{lemma:burnside})
applied to the action of $C_G(H) \times K$ on the transporter set
$N_G(H,K)$, $C_G(H) \times N_G(H,K) \times K \to N_G(H,K) \colon
(h,n,k) \to hnk$, with isotropy subgroup $C_K(H^n)$ at $n \in
N_G(H,K)$. In particular, $|\catt FG(H)| |H| = |\cat FG(H)| |Z(H)|$.

Define $\catt SG$ to be the $G$-category with objects the nonidentity
$p$-subgroups of $G$ and \m s
\begin{equation*}
  \catt SG(H,K) =
  \begin{cases}
    C_K(H) & H \leq K \\ \emptyset & H  \nleq K
  \end{cases}
\end{equation*}
with composition in $\catt SG$ induced from composition in the group
$G$.  Using that $G$ acts on $\catt SG$ by conjugation,
we may rewrite equation \eqref{eq:tildeFG1} as
\begin{equation}\label{eq:tildeFG2}
   |C_G(H)| |\catt FG(H,K)| |K| =  |(\catt SG)_{hG}(H,K)|
\end{equation}
according to \eqref{eq:ShG}.

\begin{cor}\label{cor:chiFGeqchiFtildeG}
  $\chi(\catt FG) = \chi(\cat FG)$
\end{cor}
\begin{proof}
  For any nonidentity $p$-subgroup $K$ of $G$, Theorem~\ref{thm:STLFO}
  says that
  \begin{equation*}
    \sum_{1 < H \leq K} -\mu(H) |C_K(H)| = |K| \chi(\cat FK) = |K|
  \end{equation*}
  because $\chi(\cat FK)=1$ by Proposition~\ref{prop:chiFG=1}. By
  Theorem~\ref{thm:ShG}.\eqref{prop:cowetildeF2} and equation
  \eqref{eq:tildeFG2}, $k_K = -|G|^{-1}|C_G(K)|\mu(K)$ is a co\we\ for
  $\catt FG$. But this function is also a co\we\ for $\cat FG$ by
  Theorem~\ref{thm:STLFO} and Table~\ref{tab:cat*}.
\end{proof}



\begin{lemma}[Burnside's counting lemma] \cite{MR562002} \label{lemma:burnside}
  If $X$ is a finite right $G$-set then
\begin{equation*}
   \sum_{g \in G} |X^g| = |X/G| |G|  =  \sum_{x \in X}
  |{}_xG| 
\end{equation*}
where $X^g \subset X$ is the fixed set for $g \in G$ and ${}_xG \leq
G$ is the isotropy subgroup for $x \in X$.
\end{lemma}

\subsection{Alternative \we s and co\we s}
  \label{sec:weightFG}
  We shall first reformulate the expressions for the \we s $k^\bullet$
  from Table~\ref{tab:cat*} using the \Mb\ function for $[\catn
  TG*]$ (Theorem~\ref{thm:wefrompi0C}).

  The rational number
  \begin{equation}
    \label{defn:nu}
    [\mu](H,K) = \frac{1}{|N_G(H)|} \sum_{L \in [K]} \mu(H,L)
  \end{equation}
  only depends on the conjugacy classes of $H$ and $K$\footnote{The
    function $\nu(H,K)$ is not the same as
    $|N_G(H)|^{-1}\lambda(H,K)$ where $\lambda$ is the \Mb\ function
    for the poset of $p$-subgroup classes ordered by subconjugation.}.
In particular, $\nu(K)=|N_G(K)|^{-1}\mu(K)$, where $\nu(K)$ is short
for $\nu(1,K)$.

    \begin{prop}\label{prop:nupi0TG*}
      The function $[\mu]([H],[K])$ defined by equation~\eqref{defn:nu}
      is the \Mb\ function for $[\catn TG*]$.
    \end{prop}
\begin{proof}
  We claim that $\left([\mu]([H],[K])\right)_{[H],[K] \in [\catn TG*]}$ is
  the inverse of the matrix $\left(|N_G(H,K)|\right)_{[H],[K] \in
    [\catn TG*]}$.  If $H$ and $L$ are nonidentity $p$-subgroups of
  $G$ then
  \begin{multline*}
    |N_G(H)|\sum_{[K]} [\mu]([H],[K]) |N_G(K,L)|  =
     \sum_K \mu(H,K)  |N_G(K,L)|  =
    \sum_K \mu(H,K) \sum_{g \in G} \catn SG*(K^g,L) \\  =
     \sum_{g \in G} \sum_K \mu(H,K) \catn
    SG*(K,L^{g^{-1}})  =
      \sum_{g \in G} \delta(H,L^{g^{-1}}) =
      \sum_{g \in G} \delta(H^g,L)
  \end{multline*}
  This last sum is $0$ if $H$ and $L$ are not conjugate and it is
  $|N_G(H)|$ if they are conjugate.
\end{proof}

    \begin{prop}\label{prop:wecat*}
      Weightings and \Euc s of the $p$-subgroup categories $\catn
      TG*$, $\catn LG*$, $\catn FG*$, and $\catn OG*$ are as in
      Table~\ref{tab:wecat*}.
    \end{prop}
    \begin{proof}
      In case of $\catn TG*$ this follows immediately from
      Theorem~\ref{thm:wefrompi0C} and
      Proposition~\ref{prop:nupi0TG*}. The remaining subgroup
      categories can be handled similarly. For instance, the \we\ for
      $\catn OG*$, or rather $[\catn OG*]$, is 
      \begin{equation*}
        k_{\mathcal{O}}^{[H]} = \frac{|G|}{|N_G(H)|}k_{\mathcal{O}}^H =
        \frac{|H|}{|N_G(H)|}\sum_K \mu(H,K) =
        |H| \sum_{[K]} \frac{1}{|N_G(H)|} \sum_{L \in [K]} \mu(H,L) =
        |H| \sum_{[K]} \nu([H],[K])
      \end{equation*}
      The third column of Table~\ref{tab:wecat*} is simply the sum
      $\sum_{[K]} k^{[K]}$ of the co\we s for $[\mathcal{C}]$ as in
      Theorem~\ref{thm:wefrompi0C}.
    \end{proof}

    \begin{table}[t]
      \centering
      \begin{tabular}{c||l|l|l}
        $\mathcal{C}$ & $k^{[H]}$ & $k^{[K]}$ & $\chi(\mathcal{C})$ \\ \hline
        $\catn TG*$ & 
        $\displaystyle  \sum_{[K]}
        [\mu]([H],[K])$ & 
        $\displaystyle -[\mu]([K])$ &
        $\displaystyle 
        \sum_{[K]}-[\mu]([K])$ \\
        $\catn LG*$ &
        $\displaystyle  \sum_{[K]} 
         [\mu]([H],[K]) |O^pC_G(K)|$ &
        $\displaystyle -[\mu]([K])  |O^pC_G(K)|$ &
        $\displaystyle 
        \sum_{[K]} -[\mu]([K])|O^pC_G(K)|$ \\
        $\catn FG*$ &
        $\displaystyle  \sum_{[K]} [\mu]([H],[K]) |C_G(K)|$ &
        $\displaystyle -[\mu]([K]) |C_G(K)|$ &
        $\displaystyle 
        \sum_{[K]} -[\mu]([K])|C_G(K)|$ \\
        $\catn OG*$ &
        $\displaystyle |H| \sum_{[K]}
        [\mu]([H],[K]) $ &
        $\displaystyle \sum_{[H]} |H| [\mu]([H],[K])$ &
          $\displaystyle \sum_{[H],[K]} |H| [\mu]([H],[K])$
      \end{tabular}
      \caption{Weightings, co\we s and \Euc s for categories of nonidentity
        $p$-subgroups}\label{tab:wecat*} 
    \end{table}

  \begin{rmk}\label{rmk:nu}
  Define the $\mu$-transporter from $H$ to $K$ to be the set
  \begin{equation*}
    N_G^\mu(H,K) = \{g \in G \mid \Phi(K) \leq H^g \leq K \}, \qquad 
    H,K \in \Ob{\catn SG*}
  \end{equation*}
  of group elements $g$ that conjugate $H$ into $K$ such that
  $\mu(H^g,K) \neq 0$. 
  
  The map $g \to K^{g^{-1}}$ is a bijection between
  ${N_G^\mu(H,K)}/N_G(K)$ and the set $\{ L \in [K] \mid H \leq L,
  \mu(H,L) \neq 0 \}$ of subgroups $L$ of $G$ conjugate to $K$ and
  containing $H$ with $\mu(H,L) \neq 0$ and therefore
      \begin{equation*}
        [\mu]([H],[K]) = 
        (-1)^n p^{\binom n2} \frac{ |N_G^\mu(H,K)|}{|N_G(H)| |N_G(K)|},
     \qquad H,K \in \Ob{\catn FG*}, \quad |K| = p^n |H|
      \end{equation*}
      can be computed from these transporter sets.
  \end{rmk}

  Next we note  that the values of the \we s for the $p$-subgroup
  categories of Table~\ref{tab:cat*} can be computed locally.


 Fix $H$, a nonidentity $p$-subgroup of $G$, and consider the
  projection $\catn TG*(H)= N_G(H) \to \overline{N_G(H)} = N_G(H)/H =
  \catn OG*(H)$ of the $p$-local subgroup $N_G(H)$ onto its quotient
  $N_G(H)/H$. The functor
  \begin{equation}\label{defn:CGfunctor}
    C_G \colon \big( \catn S{\catn OG*(H)}* \big)^{\mathrm{op}} \to 
    \catn {\overline{S}}{C_G(H)}{}
  \end{equation}
  takes the nonidentity $p$-subgroup $\overline{K}$ of $\catn OG*(H)$
  to the subgroup $C_G(K)$ of $C_G(H)$ where $K \leq N_G(H)$ is the
  preimage of $\overline{K} \leq N_G(H)/H$. For every $x \in C_G(H)$,
 \begin{equation}\label{defn:CG^-1}
    C_G / \langle x \rangle = 
    \{ \overline{K} \in \Ob{\catn S{\catn OG*(H)}*} \mid 
    C_{C}(K) \ni x \} =
    \catn S{\overline{C_{N_G(H)}(x)}}*
  \end{equation}
  is the preimage under $C_G$ of the subposet $\{ Y \mid \langle x
  \rangle \leq Y \leq C_G(H) \}$.  Following the bar convention of \cite[p
  18]{GLSI}, we write $\overline{C_{N_G(H)}(x)}$ for the image in
  $\catn OG*(H)$ of the centralizer in $N_G(H)$ of $x \in C_G(H)$.


  \begin{prop}\label{prop:altwes}
The \we s for $\catn
  SG*$, $\catn TG*$, and $\catn OG*$ are
  \begin{equation*}
    k_{\mathcal{S}^*}^H = 1-\chi(\catn S{\catn OG*(H)}*), \qquad
    k_{\mathcal{T}^*}^H = \frac{1-\chi(\catn S{\catn OG*(H)}*)}{|G|}, \qquad
    k_{\mathcal{O}^*}^H = 
     \frac{1-\chi(\catn S{\catn OG*(H)}*)}{|G \colon H|} 
  \end{equation*}
  and the \we\ for $\catn FG*$ is
  \begin{equation*}
     k_{\mathcal{F}^*}^H = |G|^{-1} \sum_{x \in C_G(H)} 
     \big(1- \chi(C_G/\langle x \rangle) \big)
  \end{equation*}
  \end{prop}
  \begin{proof}
    Equation \eqref{eq:altmuposet} shows that
 \begin{equation*}
  k^H_{\mathcal{S}^*} 
                  =\sum_K \mu(H,K) = \sum_{K \in [H,N_G(H)]} \mu(H,K) 
                  = 1- \chi((H,K])
                  = 1- \chi((H,N_G(H)])
                  = 1-\chi(\catn S{\catn OG*(H)}*)
\end{equation*}
as $\mu(H,K)=0$ unless $H$ is normalized by $K$
(Lemma~\ref{lemma:muH}).  (Indeed, the subposets $(H,K]$ and
$(H,N_G(H)]$ of $\catn SG*$ are homotopy equivalent \cite[Proposition
6.1]{quillen78}.)  Similarly,
\begin{multline*}
   |G|k_{\mathcal{F}^*}^H = 
    \sum_{H \leq K \leq N_G(H)} \mu(H,K) |C_G(K)| =
   \sum_{H \leq K \leq N_G(H)} \mu(H,K) |C_{N_G(H)}(K)| =
   \sum_{\overline K \leq \catn OG*(H)} 
   \mu(\overline K) |C_{N_G(H)}(K)| \\
   = |C_G(H)| - 
   \sum_{1 \lneqq \overline K \leq \catn OG*(H)} 
   -\mu(\overline K) |C_{N_G(H)}(K)|
\end{multline*}
because $C_G(K) = C_{N_G(H)}(K)$ as $C_G(K) \leq C_G(H) \leq N_G(H)$
when $H \leq K \leq N_G(H)$.  The sum that occurs in this formula for
$|G|k_{\mathcal{F}^*}^H$ is the \Euc\ \cite[Proposition
2.8]{leinster08} of the Grothendieck construction for the presheaf
$C_G$ \eqref{defn:CGfunctor}. Since the opposite of this Grothendieck
construction is the direct sum \cite[Chp 3.2]{stanley97} over $x \in
C_G(H)$ of the subposets \eqref{defn:CG^-1}
we arrive at the formula
that we wanted to prove.
\end{proof}


  Using the expressions from Proposition~\ref{prop:altwes}, the \Euc s
  of $\catn SG*$, $\catn TG*$, and $\catn OG*$ are
\begin{equation}\label{eq:chiSGchiOG}
  \chi(\catn SG*) = \sum_{H} 
  \big(1-\chi(\catn S{\catn OG*(H)}*)\big), \qquad
  \chi(\catn TG*) = 
  \sum_{[H]} \frac{1-\chi(\catn S{\catn OG*(H)}*)}{|\catn TG*(H)|},
  \qquad 
  \chi(\catn OG*) 
  =  \sum_{[H]} \frac{1-\chi(\catn S{\catn OG*(H)}*)}{|\catn OG*(H)|}
\end{equation}
The first of these equation can also be written
\begin{equation}\label{eq:combiSG}
    \sum_{H} \big(1-\chi(\catn S{\catn OG*(H)}*)+\mu(H) \big) = 0
\end{equation}
as $\chi(\catn SG*) = \sum_H -\mu(H)$ (Table~\ref{tab:cat*}).
Similarly, we obtain the alternative formula
\begin{equation*}
  \chi(\catn FG*) = |G|^{-1} \sum_H \sum_{x \in C_G(H)} 
\big(1- \chi(C_G/\langle x \rangle) \big)
\end{equation*}
for the \Euc\ of $\catn FG*$. Comparing this new formula with the one
from Table~\ref{tab:cat*} we arrive at the combinatorial identity
\begin{equation}\label{eq:combiFG}
  \sum_{H} \sum_{x \in C_G(H)} 
  \big( 1- \chi(C_G/\langle x \rangle) +\mu(H)\big) =0 
\end{equation}
where $H$ runs over the set of nonidentity $p$-subgroups of $G$.

Finally, we observe that only $p$-radical $p$-subgroups contribute to
the \we s for $\catn SG*$, $\catn TG*$, and $\catn OG*$.

\begin{defn} \label{defn:radical}  The
  $p$-subgroup $H$ of $G$ is
  \begin{itemize}
  \item $p$-radical if 
    $O_p\catn OG{}(H)=1$ \cite[Proposition 4]{bouc84a}
  \item $\catn FG{}$-radical if $O_p \widetilde{\mathcal{F}}_G(H)=1$
    \cite[Definition A.9]{blo03}
  \end{itemize}
\end{defn}

\begin{cor}\label{cor:pradical}
  The \we s for  $\catn SG*$, $\catn TG*$, and $\catn OG*$ are
  supported on the nonidentity $p$-radical subgroups of $G$.
\end{cor}
\begin{proof}
  If $O_p\catn OG*(H)>1$ then $\chi(\catn S{\catn OG*(H)}*)=1$
  \cite[Proposition 2.4]{quillen78}
  and the \we s $k^H=0$ for the categories
  $\catn SG*$, $\catn TG*$, and $\catn OG*$
  (Proposition~\eqref{prop:altwes}).
\end{proof}

Consequently,
\begin{equation*}
  \chi(\catn C{}r) = \chi(\catn C{}*), \qquad 
\mathcal{C} = \catn SG{}, \catn TG{}, \catn OG{}
  \end{equation*}
  where $\mathcal{C}^r$ is the subposet of $\mathcal{C}^*$ of
  nonidentity $p$-radical $p$-subgroups.  Bouc
  \cite[Corollaire]{bouc84b} shows the stronger result that $\catn
  SGr$ and $\catn SG*$ are homotopy equivalent posets. Th\'evenaz and
  Webb \cite[Theorem 2.3]{thevenaz_webb91} describe $\catn SGr$ when
  $G$ is simple group of Lie type in defining characteristic $p$.

We suspect that the  \we s for  $\catn SG*$, $\catn TG*$, and
$\catn OG*$ are supported  {\em
  precisely\/} on the nonidentity $p$-radical $p$-subgroups, ie that
\begin{equation}\label{conj:pradical}
  \chi(\catn S{\catn OG*(H)}*)  \neq 1 \iff O_p \catn OG*(H) = 1
\end{equation}
This would be true if the {\em strong\/} Quillen conjecture
\begin{equation}\label{conj:pradicalQuillen}
  \chi(\catn SG*) \neq 1 \iff O_pG = 1
\end{equation}
turns out to be true for all finite groups $G$.  (It is true, as used
above, that $ O_p(G) \neq 1 \Longrightarrow \chi(\cat SG)=1$ but the
problem is that $\chi(\cat SG)=1 \Longrightarrow O_p(G) \neq 1$ is
only known to hold for $p$-solvable groups with abelian \syl ps
\cite[Lemma 1.1, Theorem A]{HI88}.  The original Quillen conjecture
\cite[Conjecture 2.9]{quillen78}, that $\catn SG* \simeq \ast
\Longrightarrow O_p(G) \neq 1$, is true when $G$ is solvable.  Also,
it is known that $|G|_p$ divides $1-\chi(\cat SG)$
\cite{brown75,quillen78,thevenaz87,HIO89}.)  Explicit computations
with Magma \cite{magma} reveal that \eqref{conj:pradicalQuillen} is
true at $p=2$ for all groups of order $\leq 1500$.

We are not aware of any similar characterization of the support of the
\we\ for $\catn FG*$.

The two concepts of radical subgroups introduced in
Definition~\ref{defn:radical} are unrelated in general \cite[Appendix
A]{blo03}.  If $P$ is an abelian nonidentity $p$-group, then all
subgroups of $P$ are $\catn FG{}$-radical but only $P$ itself is
$p$-radical.  However, if $H$ is a \sfc\ $p$-subgroup of $G$
(Definition~\ref{defn:sfc}) then $O^pC_G(H)$ is a $p'$-group
(Lemma~\ref{lemma:CKHg}.\eqref{lemma:CKHg1}) and the short exact
sequence
   \begin{equation*}
     1 \to O^pC_G(H) \to \catn OGc(H) \to
     \widetilde{\mathcal{F}}_G^c(H)   \to 1
   \end{equation*} 
   can be used to verify the implication
  \begin{equation*}
    \text{$H$ is \sfc\ and $\catn FG{}$-radical} \Longrightarrow
     \text{$H$ is \sfc\ and $p$-radical}
  \end{equation*}
  The converse implication does not hold in general: Let $p=2$. The
  normal cyclic subgroup $H=O_pG$ of order $4$ in the dihedral group
  $G=D_{24}$ of order $24$ is a \sfc\ subgroup with $\catn OGc(H) =
  \Sigma_3$ and $ \widetilde{\mathcal{F}}_G^c(H)=C_2$.  Thus $H$ is
  $p$-radical but not $\catn FG{}$-radical.

\section{\Euc s of orbit categories}
\label{sec:OG}
We shall now derive a more concise expression than the ones given in
Table~\ref{tab:cat*} or Table~\ref{tab:wecat*} for the \Euc\ of $\cat
OG$.

\begin{thm}\label{thm:chiOG}
  The \Euc\ of the orbit category $\cat OG$ is
  \begin{equation*}
    \chi(\cat OG) = 
   \chi(\cat TG) + \frac{1}{|G|}\frac{p-1}{p} 
      \sum_{\text{$C \in \Ob{\cat OG}$ cyclic}} |C|
  \end{equation*}
\end{thm}
\begin{proof}
  The co\we\ $k_\bullet^{\mathcal{O}}$ from Table~\ref{tab:cat*} for
  $\cat OG$ multiplied by $|G|$ is
  \begin{equation*}
    |G|k^{\mathcal{O}}_K = \sum_{1<H}|H|\mu(H,K) 
           = -\mu(K)+ \sum_{1 \leq H}|H|\mu(H,K) =
           \begin{cases}
             -\mu(K)+\frac{p-1}{p}|K| &  \text{$K$ is cyclic} \\
             -\mu(K) & \text{$K$ is not cyclic} 
           \end{cases}
  \end{equation*}
  by Theorem~\ref{thm:STLFO} and Corollary~\ref{cor:sumKmuKL}.
  Thus the \Euc\ of  $\cat OG$ is  
  \begin{equation*}
    \chi(\cat OG) = 
   \frac{1}{|G|} \sum_{K \in \Ob{\cat OG}} -\mu(K) + 
     \frac{1}{|G|} \frac{p-1}{p} \sum_{\text{$C \in \Ob{\cat OG}$ cyclic}} 
    |C| 
  \end{equation*}
  where the first term is the \Euc\ of the transporter category $\cat
  TG$. 
 \end{proof}

  Equivalently, the \Euc\ of the orbit category $\cat OG$ is
  \begin{equation}\label{eq:chiOG2}
     \chi(\cat OG) = 
       \chi(\cat TG) +
        \frac{p-1}{p}\sum_{\text{$[C] \in [\cat FG]$ cyclic}} 
        \frac{1}{|\cat OG(C)|}
  \end{equation}
  where the sum is taken over the set of {\em conjugacy classes\/} of
  nonidentity cyclic $p$-subgroups of $G$.

\begin{cor}\label{cor:sumKmuKL}
  For any $K \in \Ob{\cat SG}$
  \begin{equation*}
    \frac{1}{|\Phi(K)|} \sum_{1 \leq H \leq K} |H| \mu(H,K)=
    \begin{cases}
      p-1 & \text{$K$ is cyclic} \\ 0  & \text{$K$ is not cyclic}
    \end{cases}
  \end{equation*}
\end{cor}
\begin{proof}
  Suppose that the Frattini quotient $K/\Phi(K)$ is elementary abelian
  of order $p^{n}$ for some $n > 0$. Recall that $n=1$, $K/\Phi(K)$ is
  cyclic, if and only if $K$ is cyclic \cite[Chp 5, Corollary
  1.2]{gorenstein68}.  The sum of this corollary,
\begin{equation*}
   \sum_{H \colon \Phi(K) \leq H \leq K}
  |H/\Phi(K)| \mu(K/H) =
  \sum_{d=0}^{n} (-1)^{n-d}\binom{[n]}{[d]}
  p^{\binom{n-d}2} p^d 
  =   \sum_{d=0}^{n} (-1)^{d} 
  \binom{[n]}{[d]} p^{\binom{d}2} p^{n-d},
\end{equation*}
is evaluated in Lemma~\ref{lemma:sum}. It is nontrivial only if $n=1$
where it has value $p-1$.
\end{proof}

The Gaussian $p$-binomial coefficient
\begin{equation*}
\binom{[n]}{[d]} = \frac{\prod_{j=1}^d\big(p^{n}-p^{j-1}\big)}
              {\prod_{j=1}^d\big(p^{d}-p^{j-1}\big)} =
\frac{\prod_{j=1}^d\big(p^{n+1-j}-1\big)}
              {\prod_{j=1}^d\big(p^{j}-1\big)} 
\end{equation*}
counts the number of $d$-dimensional subspaces of the $n$-dimensional
$\F_p$-vector space $\F_p^{n}$ \cite[1.3.18]{stanley97}.

\begin{lemma}\label{lemma:sum}
  For any $n \geq 1$,
  \begin{equation*}
     \sum_{d=0}^{n} (-1)^{d} 
  \binom{[n]}{[d]} p^{\binom{d}2} p^{n-d} =
  \begin{cases}
    p-1 & n=1 \\ 0 & n>1
  \end{cases}
  \end{equation*}
\end{lemma}
\begin{proof}
Note first the formulas \cite[p 26]{stanley97}
\begin{equation*}
  \binom{[n]}{[d]} =  \binom{[n-1]}{[d]}+p^{n-d}\binom{[n-1]}{[d-1]}, \qquad 
  \binom{[n]}{[0]} = 1 =  \binom{[n]}{[n]}, \qquad
   \binom{[2]}{[1]} = 1+p,
\end{equation*}
for the Gaussian $p$-binomial coefficients. 

For $n=1$ and $n=2$, the sums we are evaluating are the polynomials
\begin{equation*}
  \binom{[1]}{[0]}p -  \binom{[1]}{[1]} = p-1, \qquad
   \binom{[2]}{[0]}p^2 -  \binom{[2]}{[1]}p +  \binom{[2]}{[2]}p =
   p^2-(1+p)p+p = 0
\end{equation*}
For $n>2$ the sum has the value
\begin{align*}
  \sum_{d=0}^{n} (-1)^{d} 
  \binom{[n]}{[d]} p^{\binom{d}2} p^{n-d} 
  &= p^n +  \sum_{d=1}^{n-1} (-1)^{d} 
  \binom{[n]}{[d]} p^{\binom{d}2} p^{n-d} + (-1)^np^{\binom n2} \\
  &= p^n +  \sum_{d=1}^{n-1} (-1)^{d} 
  \left( \binom{[n-1]}{[d]} + p^{n-d} \binom{[n-1]}{[d-1]} \right) 
    p^{\binom{d}2} p^{n-d} + (-1)^np^{\binom n2} \\
   &=\left(p^n +  \sum_{d=1}^{n-1} (-1)^{d}  \binom{[n-1]}{[d]}
     p^{\binom{d}2} p^{n-d}\right) +
   \left(  \sum_{d=1}^{n-1} (-1)^{d}  \binom{[n-1]}{[d-1]}
     p^{\binom{d}2} p^{2(n-d)} + (-1)^np^{\binom n2} \right) \\
   &= \sum_{d=0}^{n-1} (-1)^{d}  \binom{[n-1]}{[d]}
     p^{\binom{d}2} p^{n-d} +
      \sum_{d=1}^{n} (-1)^{d}  \binom{[n-1]}{[d-1]}
     p^{\binom{d}2} p^{2(n-d)} 
\end{align*}
The first term is
\begin{equation*}
  \sum_{d=0}^{n-1} (-1)^{d}  \binom{[n-1]}{[d]}
     p^{\binom{d}2} p^{n-d} 
   = p \sum_{d=0}^{n-1} (-1)^{d}  \binom{[n-1]}{[d]}
     p^{\binom{d}2} p^{n-1-d}   
\end{equation*}
and the second term is 
\begin{align*}
   \sum_{d=1}^{n} (-1)^{d}  \binom{[n-1]}{[d-1]}
     p^{\binom{d}2} p^{2(n-d)}
   &= -   \sum_{d=0}^{n-1} (-1)^{d}  \binom{[n-1]}{[d]}
     p^{\binom{d+1}2} p^{2(n-d-1)} 
   =  -   \sum_{d=0}^{n-1} (-1)^{d}  \binom{[n-1]}{[d]}
     p^{\binom{d}2}p^d p^{2(n-d-1)} \\
    &=  -   \sum_{d=0}^{n-1} (-1)^{d}  \binom{[n-1]}{[d]}
     p^{\binom{d}2} p^{2n-d-2} 
      =  - p^{n-1}  \sum_{d=0}^{n-1} (-1)^{d}  \binom{[n-1]}{[d]}
     p^{\binom{d}2} p^{n-1-d}
\end{align*}
We have now proved the recursive relation
\begin{equation*}
    \sum_{d=0}^{n} (-1)^{d} 
  \binom{[n]}{[d]} p^{\binom{d}2} p^{n-d} =
  p(1-p^{n-2})  \sum_{d=0}^{n-1} (-1)^{d} 
  \binom{[n-1]}{[d]} p^{\binom{d}2} p^{n-1-d}
\end{equation*}
for $n>2$. Since the sum equals $0$ for $n=2$, it equals $0$ for all
$n \geq 2$. 
\end{proof}    


Equating the two expressions for $\chi(\catn OG*)$ from
Theorem~\ref{thm:chiOG} and \eqref{eq:chiSGchiOG} we arrive at the
combinatorial identity
\begin{equation}
  \label{eq:combiOG}
  \sum_H \big( |H| - \chi(\catn S{\catn OG*(H)}*)|H|+\mu(H) \big) =
  \frac{p-1}{p}\sum_{C} |C|
\end{equation}
where $H$ runs over the set of nonidentity $p$-subgroups of $G$ and
$C$ over the set of nonidentity cyclic $p$-subgroups of $G$.

\begin{cor}\label{cor:chiOGinZ}
  $|G|_{p'}\chi(\catn OG*)$ is an integer.
\end{cor}
\begin{proof}
  In fact, all values of
  \begin{equation*}
    |G|_{p'} k^{[H]}_{\mathcal{O}} =
    |G|_{p'} \frac{1-\chi(\catn S{N_G(H)/H}*)}{|N_G(H) \colon H|} =
    \frac{|G|_{p'}}{|N_G(H) \colon H|_{p'}}
    \frac{1-\chi(\catn S{N_G(H)/H}*)}{|N_G(H) \colon H|_p} 
  \end{equation*}
  are integers because $|G|_p$, and hence also $|N_G(H) \colon H|_p$,
  divides $1-\chi(\catn S{N_G(H)/H}*)$ \cite[Corollary 2]{brown75}.
\end{proof}

\begin{rmk}[\Euc\ of $\mathcal{O}_G$]\label{rmk:chiOG}
 The function
  \begin{equation*}
    k^H =
    \begin{cases}
      |G|^{-1}-\chi(\cat TG) & H=1 \\ |H| \sum_{K>1} \mu(H,K) & H>1
    \end{cases}
  \end{equation*}
  is a weighting for $\mathcal{O}_G$, and 
  \begin{equation*}
    \chi(\mathcal{O}_G) =  |G|^{-1}-\chi(\cat TG) + \chi(\cat OG) =
        |G|^{-1} +  
        \frac{p-1}{p}\sum_{\text{$[C] \in [\cat FG]$ cyclic}}
         \frac{1}{|\cat OG(C)|}
  \end{equation*}
  is the \Euc\ of the orbit category $\mathcal{O}_G$ of $G$.
\end{rmk}




\section{The range of $\chi(\cat FG)$}
\label{sec:rangechiFG}
We shall first identify a class of finite groups $G$ for which 
$\chi(\catn FG*)=1$.

\begin{prop}\label{prop:chiFG=1}
  If $G$ contains a central nonidentity $p$-subgroup then
  $\chi(\cat FG)=1$ and $\chi(\cat LG) = |G  \colon O^pG|^{-1}$.
\end{prop}
\begin{proof}
  We concentrate on $\cat FG$ and leave the similar case of $\cat LG$
  to the reader.  Let $Z$ be a nonidentity central $p$-subgroup of $G$
  and $Z^+$ the full subcategory of $\cat FG$ generated by
  $p$-subgroups containing $Z^+$. $Z^+$ is a deformation retract of
  $\cat FG$ in the sense that there are functors
  \begin{equation*}
    \xymatrix@1{Z^+ \ar@<.5ex>[r]^R & {\cat FG}
    \ar@<.5ex>[l]^L},  \qquad 1_{Z^+} = LR,
   \quad 1_{\cat FG} \Longrightarrow RL,
  \end{equation*}
  where $R$ is the inclusion functor and $L$ is the functor that takes
  $Q \leq G$ to $LQ=QZ$ and the $\cat FG$-\m\ \func{c_g}PQ 
  to \func{c_g}{LP}{LQ} (where $c_g \colon x \mapsto x^g = g^{-1}xg$
  is conjugation by $g \in G$). If $P$ and $Q \geq Z$ are nonidentity
  $p$-subgroups of $G$ then
  \begin{equation*}
    Z^+(LP,Q) = \cat FG(PZ,Q) = \qt{C_G(PZ)}{N_G(PZ,Q)} =
    \qt{C_G(P)}{N_G(P,Q)} = \cat FG(P,Q) = \cat FG(P,RQ)
  \end{equation*}
  showing that $L$ and $R$ are adjoint functors with $L \dashv R$.  By
  Lemma~\ref{lemma:adjointcats}, the EI-categories $Z^+$ and $\cat FG$
  have the same \Euc s and, by Example~\ref{exmp:finalelement},
  $\chi(Z^+)=1$ as $Z^+$ has initial object $Z$.
\end{proof}

The converse of Proposition~\ref{prop:chiFG=1} is not true as
$\chi(\cat FG)=1$ and $Z(G)=1$ for $G=\Sigma_3$ and $p=2$.

The \Euc\  $\chi(\cat FG) = |G|^{-1} \sum k_K = -|G|^{-1} \sum_K
\mu(K) |C_G(K)|$ is a rational number such that $|G|\chi(\catn FG*)$
is an integer. We now improve this observation.

  \begin{cor}\label{cor:mchiGinZ}
    $|G|_{p'}\chi(\cat FG)$ is an integer.
  \end{cor}
  \begin{proof}
    In fact, all values of 
    \begin{equation*}
       |G|_{p'}k^{\mathcal{F}}_{[K]} = 
      |G|_{p'}\frac{\mu(K)}{|\cat FG(K)|} = 
      \frac{\mu(K)}{|\cat FG(K)|_p} \cdot 
      \frac{|G|_{p'}}{|\cat FG(K)|_{p'}}, \qquad K \in \Ob{\catn SGa},
      \quad |K|=p^n,
    \end{equation*}
    are integers because $|\cat FG(K)|_p$ divides $|\Aut{}{K}|_p =
    p^{\binom n2} = \mu(K)$ as $|\cat FG(K)|$ divides $|\Aut{}{K}|$,
    and $|\cat FG(K)|_{p'}$ divides $|G|_{p'}$ as $|\cat FG(K)| =
    |N_G(K) \colon C_G(K)|$ divides $|G|$. (Remember that $\mu(K)=0$
    unless $ K \in \Ob{\catn SGa}$ is elementary abelian of order
    $|K|=p^n$, $n>0$.)
  \end{proof}

We now show that the computation of the \Euc\ of the
Frobenius category can be reduced to the computation of \Euc s of
posets.

\begin{prop}\label{prop:chiSCG}
 $\displaystyle
  \chi(\cat FG) = |G|^{-1} \sum_{x \in G} \chi(\cat S{C_G(x)})$
\end{prop}
\begin{proof}
  Recall that $\overline{\mathcal{S}}_G$ denotes the poset of all
  subgroups of $G$.  Note that for any subgroup $K$ of $G$
  \begin{equation*}
    |C_G(K)| =
    |\{ x \in G \mid x \in C_G(K) \}| =
    |\{ x \in G \mid C_G(x) \geq  K \}| =
    \sum_{x \in G} \overline{\mathcal{S}}_G(K,C_G(x))
  \end{equation*}
  and therefore 
  \begin{multline*}
    |G| \chi(\catn FG*)
    = \sum_{K \in \Ob{\mathcal{S}_G^*}} -\mu(K)|C_G(K)| 
    = \sum_{K \in \Ob{\mathcal{S}_G^*}}\sum_{x \in G}
    -\mu(K)\overline{\mathcal{S}}_G(K,C_G(x)) \\
    =\sum_{x \in G}\sum_{K \in \Ob{\mathcal{S}_G^*}}
    -\mu(K)\overline{\mathcal{S}}_G(K,C_G(x))  
    = \sum_{x \in G}\sum_{K \in \Ob{\mathcal{S}_{C_G(x)}^*}} -\mu(K)
    = \sum_{x \in G} \chi(\mathcal{S}_{C_G(x)}^*)
  \end{multline*}
  where the final equality uses the formula from Table~\ref{tab:cat*}
  for $\chi(\cat SG)$.
\end{proof}

\begin{cor}\label{cor:normalsyl}
  Suppose that $G$ has a normal \syl p, $P$, (so that $G = P \rtimes
  G/P$) and let $k_{\mathcal{F}}^H$, $H \in \Ob{\catn SG*}$, be the \we\ for $\catn
  FG*$ from Table~\ref{tab:cat*}.
  \begin{enumerate}
  \item $\displaystyle |G|k_{\mathcal{F}}^H = |\{x \in G \mid
      C_P(x)=H \}|$
    \label{cor:normalsyl1}
  \item $k_{\mathcal{F}}^H \geq 0$ and $k_{\mathcal{F}}^H > 0$ if and
    only if $H=C_P(x)$ for some $x \in G$ \label{cor:normalsyl2}
  \item $\displaystyle 
    |G|\chi(\cat FG) =|\{ x \in G \mid C_P(x)>1 \}|$
    \label{cor:normalsyl3} 
  \end{enumerate}
\end{cor}
\begin{proof}
  For any nonidentity $p$-subgroup $K$ and any element $x \in G$,
  since $K \leq P$,
  \begin{equation*}
    x \in C_G(K) \iff K \leq C_G(K) \iff K \leq P \cap C_G(x) = C_P(x)
  \end{equation*}
  so that
  \begin{equation*}
    |G| k_{\mathcal{F}}^H = \sum_K \mu(H,K) |C_G(K)| = 
    \sum_{x \in G}\sum_K \mu(H,K) \catn SG*(K,C_P(x)) =
     \sum_{x \in G} \delta(H,C_p(x))
  \end{equation*}
  by the \Mb\ inversion formula \eqref{eq:zetamu}. This proves
  \eqref{cor:normalsyl1} which immediately implies \eqref{cor:normalsyl2}
  and \eqref{cor:normalsyl3}.
\end{proof}

\begin{exmp}\label{exmp:weonnoncentric}
  Let $p=2$ and $G = P \rtimes C_3$ where the cyclic group $C_3$
  cyclically permutes the three factors of $P=C_2^3$. Then $\chi(\catn
  FG*)=1$ by Proposition~\ref{prop:chiFG=1}; indeed,
  $k_{\mathcal{F}}^{Z(G)} = 2/3$, $k_{\mathcal{F}}^P = 1/3$, and
  $k_{\mathcal{F}}^H=0$ for all other nonidentity $2$-subgroups $H
  \leq G$ by Corollary~\ref{cor:normalsyl}.
\end{exmp}

\begin{cor}\label{cor:abelP}
  Suppose that $G$ has an abelian \syl p, $P$. Then
  \begin{equation*}
    \chi(\cat FG) = 
    \frac{|\{ \varphi \in \mathcal{F}_G(P) \mid C_P(\varphi)>1 \}}
    {|\mathcal{F}_G(P)| }
  \end{equation*}
\end{cor}
\begin{proof}
  When $P$ is abelian, $\catn FG*(P)$ has order prime to $p$ and 
  \begin{equation*}
    \mathcal{F}_G = \mathcal{F}_{N_G(P)} = \mathcal{F}_{P \rtimes
      \mathcal{F}_G(P)}  
  \end{equation*}
  where the first identity is Burnside's Fusion Theorem \cite[Lemma
  16.9]{GLSII} which says that $N_G(P)$ controls $p$-fusion in $G$.
  For the second equality, observe that all \m s in the Frobenius
  category of $N_G(P)$ extend to auto\m s of $P$.  Now apply
  Corollary~\ref{cor:normalsyl}.\eqref{cor:normalsyl3} to $P \rtimes
  \mathcal{F}_G(P)$.
\end{proof}

For instance, let $D_{pn}$ be the dihedral group of order $2pn$, $n
\geq 1$, $A_p$ the alternating group of order $p>2$, and
$\SL{}2{\F_q}$ the special linear group where $q$ is a
power of $p$.  Then
  \begin{equation*}
    \chi(\cat F{D_{pn}}) = \frac{1}{(2,p-1)}, \qquad
    \chi(\cat F{A_p}) = \frac{2}{p-1}, \qquad
    \chi(\cat F{\SL{}2{\F_q}}) = \frac{(2,q-1)}{q-1}
  \end{equation*}
  The computer-generated Table~\ref{tab:FAn} displays \Euc s of
  Frobenius categories at $p=2$ of small alternating groups.  (The
  Frobenius categories for $A_{2n}$ and $A_{2n+1}$ at $p=2$ are
  equivalent.)  We do not know if the sequence $\chi(\cat F{A_n})$
  converges.


  \begin{exmp}\label{exmp:chi>1}
    The group $H=(C_3 \times C_3) \rtimes C_2$, where $C_2$ swaps the
    two copies of $C_3$, has an irreducible $4$-dimensional
    representation over $\F_2$. Let $G=C_2^4 \rtimes H$ be the
    associated semi-direct product. Then $|G|=288$ and $\chi(\cat
    FG)=10/9$ at $p=2$.
  \end{exmp}


  \begin{table}[t]
    \centering
    \begin{tabular}{c||c|c|c|c|c|c|c|c|c|c|c} \hline
      $n$ & $4$ & $5$ & $6$ & $7$ & $8$ & $9$ & $10$ & $11$ & $12$ & 
      $13$ & $14$ \\
      \hline \hline
      $\chi(\cat S{A_n})$ & $1$ & $5$ & $-15$ & $-175$ & $65$ & $5121$ &
      $15105$ & $55935$ & $-288255$ & $1626625$ & $23664641$ \\ \hline
      $\chi(\cat L{A_n})$ 
      & \multicolumn{2}{|c}{$1/12$} 
      & \multicolumn{2}{|c}{$-1/24$} 
      & \multicolumn{2}{|c}{$1081/2016$} 
      & \multicolumn{2}{|c}{$971/6720$} 
      & \multicolumn{2}{|c|}{$90035/145152$} 
      & $406699/1451520$ \\ \hline
      $\chi(\cat F{A_n})$ 
      & \multicolumn{2}{|c}{$1/3$}  
      & \multicolumn{2}{|c}{$1/3$}  
      & \multicolumn{2}{|c}{$41/63$}  
      & \multicolumn{2}{|c}{$18/35$} 
      & \multicolumn{2}{|c|}{$389/567$} 
      & $233/405$ \\ \hline
      $\chi(\cat O{A_n})$ & $1/3$ & $1/3$ & $1/3$ & $2/9$ & $13/63$ & 
      $44/315$ & $178/945$ & $46/315$ & $397/2835$ & $160/1701$ &
      $2101/42525$ \\ \hline
    \end{tabular}
    \caption{\Euc s at $p=2$ of  nonidentity $p$-subgroup posets and
      categories of alternating groups computed by Magma \cite{magma}}
    \label{tab:FAn}
  \end{table}

  In all the examples we have checked (and also if $G$ has a normal or
  abelian \syl p as in Corollary~\ref{cor:normalsyl} and
  Corollary~\ref{cor:abelP}) $\chi(\cat FG)$ is positive when $p$
  divides the order of $G$.  Example~\ref{exmp:chi>1} shows that
  $\chi(\cat FG)$ can be greater than $1$. Prompted by these
  observations, we would like to pose two questions:
  \begin{itemize}
  \item Is $\chi(\cat FG)$ always positive when $p$ divides the order of $G$?
  \item Can  $\chi(\cat FG)$ get arbitrarily large?
  \end{itemize}

\section{Product formulas}
\label{sec:prodform}

We present product formulas for the \Euc s of the subgroup poset
$\catn S{G_1 \times G_2}*$ and the Frobenius category $\catn F{G_1
  \times G_2}*$ for the product of two finite groups $G_1$ and $G_2$. 


The formula of Theorem~\ref{thm:STLFO} for $\chi(\cat SG)$ may be
written in the alternative form of \eqref{eq:altmuposet} as
  \begin{equation*}
    1-\chi(\mathcal{S}_G^*) =
    \sum_{H \in \Ob{\mathcal{S}_G}} \mu(H)
    \end{equation*}
    with summation over {\em all\/} $p$-subgroups of $G$.  We shall
    use this expression to derive a formula for the Euler
    characteristic of the subgroup poset of a direct product of
    groups.

\begin{thm}\label{thm:chiSprod}
  Let $G_1,\ldots,G_n$ be finite groups. Then
  \begin{equation*}
    1-\chi(\mathcal{S}_{\prod_{i=1}^n G_i}^*) = 
    \prod_{i=1}^n 1-\chi(\mathcal{S}_{G_i}^*)
  \end{equation*}
\end{thm}
\begin{proof}
  By induction over $n$ it is enough to prove the formula for a
  product of two groups, $G_1$ and $G_2$. It is then equivalent to
  \begin{equation*}
    \sum_{H \in \Ob{\mathcal{S}_{G_1 \times G_2}}} \mu(H) =
     \sum_{H_1 \in \Ob{\mathcal{S}_{G_1}}} \mu(H_1) \;\;\cdot
       \sum_{H_2 \in \Ob{\mathcal{S}_{G_2}}} \mu(H_2) 
  \end{equation*}
  Let $\pi_1 \colon G_1 \times G_2 \to G_1$ and $\pi_2 \colon G_1
  \times G_2 \to G_2$ be the projections.  The product poset
  $\mathcal{S}_{G_1} \times \mathcal{S}_{G_2}$ \cite[Chp
  3.2]{stanley97} is a deformation retract of $\mathcal{S}_{G_1 \times
    G_2}$ in the sense that there are poset \m s
  \begin{equation*}
    \xymatrix@1{
      {\mathcal{S}_{G_1 \times G_2}} \ar@<0.5ex>[r]^-L &
      {\mathcal{S}_{G_1} \times \mathcal{S}_{G_2}} \ar@<0.5ex>[l]^-R},
    \qquad 1_{\mathcal{S}_{G_1} \times \mathcal{S}_{G_1}}= LR, \quad
    1_{\mathcal{S}_{G_1 \times G_2}} \Longrightarrow RL,
  \end{equation*}
  where $LH = (\pi_1(H),\pi_2(H))$, $H \leq G_1 \times G_2$, 
  $R(H_1,H_2)=H_1 \times H_2$, $H_1 \leq G_1$, $H_2 \leq G_2$, and
  \begin{equation*}
    H \leq R(H_1,H_2) \iff LH \leq (H_1,H_2)
  \end{equation*}
  In this situation
  \begin{equation*}
    \sum_{H \in \Ob{\mathcal{S}_{G_1 \times G_2}} 
      \colon LH=(H_1,H_2)} \mu(H) = 
      \mu(H_1) \mu(H_2) 
  \end{equation*}
  for all $p$-groups $H_1 \leq G_1$, $H_2 \leq G_2$, and $H \leq G_1
  \times G_2$ by \cite[Proposition 4.4]{leinster08} and the formula
  \cite[Proposition 3.8.2]{stanley97}
  \begin{equation*}
    \mu_{\mathcal{S}_{G_1} \times
  \mathcal{S}_{G_2}}((1,1),(H_1,H_2)) = \mu(H_1) \mu(H_2)
  \end{equation*}
  for the \Mb\ function $\mu_{\mathcal{S}_{G_1} \times
    \mathcal{S}_{G_2}}$ of the product poset $\mathcal{S}_{G_1} \times
  \mathcal{S}_{G_2}$. The theorem now easily follows.
\end{proof}
  
At $p=2$, $1-\chi(\cat S{\Sigma_6}) = 16$ and $1-\chi(\cat
S{\Sigma_6/A_6}) = 1-\chi(\cat F{C_2}) = 0$, where $\Sigma_6 = A_6
\rtimes C_2$ is the permutation group and $A_6$ the alternating group.
This example shows that Theorem~\ref{thm:chiSprod} does not generalize
to semi-direct products.



Theorem~\ref{thm:chiSprod} also follows from work of Quillen.
According to \cite[Proposition 2.6]{quillen78}, $\catn S{G_1 \times
  G_2}a$ is homotopy equivalent to the join $\catn S{G_1}a \ast \catn
S{G_2}a$ and therefore
\begin{equation*}
  1-\chi(\catn S{G_1 \times G_2}a) = 1-\chi(\catn S{G_1}a \ast \catn S{G_2}a) 
  = \big(1-\chi(\catn S{G_1}a)\big)  \big(1-\chi(\catn S{G_2}a)\big)
\end{equation*}
as  $1-\chi(X \ast Y) = \big(1-\chi(X)\big)\big(1-\chi(Y)\big)$ for
any two finite abstract simplicial complexes, $X$ and $Y$.


The formula of Theorem~\ref{thm:STLFO} for $\chi(\cat FG)$ may be
rewritten as
\begin{equation*}
  1-\chi(\cat FG) =  \frac{1}{|G|}
  \sum_{H \in \Ob{\mathcal{S}_G}} \mu(H)|C_G(H)| 
\end{equation*}
with summation over {\em all\/} $p$-subgroups $H$ of $G$. We shall use
this expression to derive a formula for the Euler characteristic of
the Frobenius category of a direct product of groups.

\begin{thm}\label{thm:chiFprod}
  Let $G_1,\ldots,G_n$ be finite groups. Then
  \begin{equation*}
    1-\chi(\mathcal{F}_{\prod_{i=1}^n G_i}^*) = 
    \prod_{i=1}^n 1-\chi(\mathcal{F}_{G_i}^*)
  \end{equation*}
\end{thm}
\begin{proof}
    By induction over $n$ it is enough to prove the formula for a
  product of two groups, $G_1$ and $G_2$. It is then equivalent to
  \begin{equation*}
    \sum_{H \in \Ob{\mathcal{S}_{G_1 \times G_2}}} \mu(H) 
    |C_{G_1 \times G_2}(H)| =
     \sum_{H_1 \in \Ob{\mathcal{S}_{G_1}}} \mu(H_1) |C_{G_1}(H_1)| \;\;\cdot
       \sum_{H_2 \in \Ob{\mathcal{S}_{G_2}}} \mu(H_2) |C_{G_2}(H_2)|  
  \end{equation*}
  But this follows as in the proof of Theorem~\ref{thm:chiSprod}
  because $C_{G_1 \times G_2}(H) = C_{G_1}(H_1) \times C_{G_2}(H_2)$
  when $H \leq G_1 \times G_2$ and $H_1=\pi_1(H)$, $H_2=\pi_2(H)$ are
  the projections of $H$.
\end{proof}

At $p=2$, $1-\chi(\cat F{\Sigma_4}) = 1/3$ and $1-\chi(\cat
F{\Sigma_4/A_4}) = 1-\chi(\cat F{C_2}) = 0$, where $\Sigma_4$ is the
permutation group and $A_4$ the alternating group. This example shows
that Theorem~\ref{thm:chiFprod} does not generalize to semi-direct
products. (Note also that $O_2(\Sigma_4) > 1$ and $\chi(\cat
F{\Sigma_4}) \neq 1$ in contrast to \cite[Lemma 1.1]{HI88} according
to which $\chi(\cat SG)=1$ whenever $O_p(G)>1$.)

\section{Variations on Frobenius categories}
\label{sec:varFrob}


Recall that $\overline{\mathcal{S}}_G$ is the poset of all subgroups of $G$.  
For any two subgroups, $H$ and $K$, of $G$, 
\begin{equation}\label{eq:zetaG}
  \overline{\mathcal{S}}_G(H,K) =
  \begin{cases}
    1 & \text{if $H \leq K$} \\ 0 & \text{if $H \nleq K$}
  \end{cases}
\end{equation}
Writing $[H]$ for the $G$-conjugacy class of $H$, let
\begin{equation*}
  \overline{\mathcal{S}}_G([H],K) =  \sum_{H \in [H]} \overline{\mathcal{S}}_G(H,K)
\end{equation*}
denote the number of subgroups of $K$ that are $G$-conjugate to $H$. In
particular, $\overline{\mathcal{S}}_G([H],G) = |[H]| = |G \colon N_G(H)|$ is the number
of conjugates of $H$ in $G$.

We next formulate an alternative expression for the \Euc\ of a
Frobenius category. Let $P$ be a subgroup of $G$ of index prime to
$p$, for instance, a \syl p of $G$.  Write $P \cap \cat FG$ for the
full subcategory of $\cat FG$ generated by all nonidentity
$p$-subgroups of $P$. Then $P \cap \cat FG$ and $\cat FG$ are
equivalent so they have identical \Euc s
(Lemma~\ref{lemma:adjointcats}).

\begin{cor}\label{cor:coweight}
  The function 
  \begin{equation*}
    k_K = 
    -\frac{\mu(K)}{|\cat FG(K,P)|}, \qquad 
    K \in \Ob{P \cap \cat FG}, 
  \end{equation*}
  is a co\we\  for $P \cap \cat FG$ and the \Euc\ of $P \cap
  \cat FG$ is
  \begin{equation*}
    \chi(P \cap \cat FG) = \sum_{K \in \Ob{P \cap \cat FG}} 
    -\frac{\mu(K)}{|\cat FG(K,P)|} = \chi(\cat FG)
  \end{equation*}
  with summation over the nonidentity elementary abelian
  $p$-subgroups $K$ of $P$.
\end{cor}
\begin{proof}
As $k_K=-{\mu(K)}/{|G \colon C_G(K)|}$ is a co\we\ for $\cat FG$
(Table~\ref{tab:cat*}), the function
\begin{equation*}
  k_K \frac{\overline{\mathcal{S}}_G([K],G)}{\overline{\mathcal{S}}_G([K],P)} 
  = -\frac{\mu(K)}{|\cat FG(K,P)|} 
  = -\frac{1}{|[K]|} \frac{\mu(K)}{|\cat FG(K)|} 
\end{equation*}
is a co\we\ for $P \cap \cat FG$. The identities
$\overline{\mathcal{S}}_G([K],G)=|G \colon N_G(K)|$ and
$\overline{\mathcal{S}}_G([K],P) = {|\cat FG(K,P)|}/{|\cat FG(K)|}$ go
into the above equality sign.
\end{proof}

Let $P$ be a finite $p$-group and $\mathcal{F}$ an abstract Frobenius
category over $P$ \cite[Chp 2]{puig09} \cite{blo03}. The objects of
$\mathcal{F}$ are the subgroups of $P$. Define $\mathcal{F}^*$ to be
the full subcategory of $\mathcal{F}$ generated by all {\em
  nonidentity\/} subgroups of $P$.
  \begin{thm}\label{thm:abscoweight}
    The function
    \begin{equation*}
      k_K = \frac{-\mu(K)}{|\mathcal{F}^*(K,P)|}, \qquad 
      K \in \Ob{\mathcal{F}^*},
    \end{equation*}
    is a co\we\  for $\mathcal{F}^*$ and the \Euc\ of
    $\mathcal{F}^*$ is
    \begin{equation*}
      \chi(\mathcal{F}^*) = \sum_{[K] \in [\mathcal{F}^*]} 
                            \frac{-\mu(K)}{|\mathcal{F}^*(K)|}
    \end{equation*}
    The \Euc\ $\chi(\catn F{}*) \in \Z_{(p)}$ is a $p$-local integer.
  \end{thm}
  \begin{proof}
    The Divisibility Axiom \cite[2.3.1]{puig09} implies that
  \begin{equation*}
    |\mathcal{F^*}(H,K)| = |\mathcal{F^*}(H)| \catn SP*([H],K)
  \end{equation*}
  where $[H] \subset \Ob{\mathcal{F}} = \Ob{\mathcal{F}_P}$ is the set
  of $\mathcal{F}$-objects $\mathcal{F}$-isomorphic to $H$ and
  $\catn SP*([H],K) = \sum_{H \in [H]} \catn SP*(H,K)$ is the number of
  $\mathcal{F}$-objects $\mathcal{F}$-isomorphic to $H$ and contained
  in $K$. In particular, $\catn SP*([H],P) = |[H]|$.  For any object $K$
  of $\Ob{\mathcal{F}^*}$
  \begin{align*}
    &\sum_{H \in \Ob{\mathcal{F}^*}}
    -\frac{\mu(H)}{|\mathcal{F}^*(H,P)|}|\mathcal{F}^*(H,K)|
    =  \sum_{H \in \Ob{\mathcal{F}^*}} -\mu(H)
    \frac{\catn SP*([H],K)}{|[H]|} 
    = \sum_{[H] \in [\mathcal{F}^*]} \sum_{H \in [H]}  -\mu(H)
    \frac{\catn SP*([H],K)}{|[H]|}  \\
    &=  \sum_{[H] \in [\mathcal{F}^*]} -\mu(H)\catn SP*([H],K)
    = \sum_{H \in \Ob{\mathcal{F}^*}} -\mu(H) \catn SP*(H,K)
    = \sum_{H \in \Ob{\cat SK}} -\mu(H) = \chi(\cat SK) = 1
  \end{align*}
  This shows that 
  \begin{equation*}
    k_K = \frac{-\mu(K)}{|\mathcal{F}^*(K,P)|}
        = \frac{-1}{|[K]|}\frac{\mu(K)}{|\mathcal{F}^*(K)|}, 
    \qquad K \in \Ob{\mathcal{F}^*}, 
  \end{equation*}
  is a co\we\ for $\mathcal{F}^*$. The formula for the \Euc\ follows.

  Since $|\mathcal{F}^*(K)|_p$ divides $|\Aut {}K|_p = \mu(K)$ when
  $K$ is elementary abelian (see proof of
  Corollary~\ref{cor:mchiGinZ}), the \Euc\ of $\mathcal{F}^*$ is a
  $p$-local integer.
  \end{proof}

\section{Self-centralizing subgroups}
\label{sec:scFG}

This section deals with the $p$-subgroup categories generated by the
\sfc\ subgroups. We mention here some facts to justify our interest in
these subcategories of \sfc\ subgroups:
\begin{itemize}
\item The centric linking category $\catn LGc$ is a
complete algebraic invariant of the $p$-completed classifying space of
$G$ \cite[Theorem A]{blo1}
\item The Frobenius category $\catn FG*$ is
completely determined by its centric subcategory $\catn FGc$ \cite[Chp
4--5]{puig09}
\item All \m s in  the category $\cattc FG$
are epi\m s \cite[Corollary
4.9]{puig09}
\item All \m s in  the category $\cattc FG$ have  unique maximal
  extensions \cite{puig09}
\end{itemize} 
Now follow the definition and a few standard properties of \sfc\
$p$-subgroups.

\begin{defn}\cite[4.8.1]{puig09} \cite[Definition A.3]{blo03} \label{defn:sfc}
  The $p$-subgroup $H$ of $G$ is \sfc\ if $C_H(H) \to C_G(H)$ is a
  $p$-Sylow inclusion.
\end{defn}

\begin{lemma}\label{lemma:CKHg}\cite[Chp 4]{puig09} \cite[Appendix A]{blo03}
  Let $H$ be a $p$-subgroup of $G$ and let $P$ be a \syl p of
  $G$.
  \begin{enumerate}
  \item $H$ is \sfc\ if and only if $C_G(H) \cong Z(H) \times O^pC_G(H)$
  and $O^pC_G(H)$ is a $p'$-group. \label{lemma:CKHg1}
  \item $H$ is \sfc\ if and only if $C_P(H^g) \leq
    H^g$ for every $g \in N_G(H,P)$. \label{lemma:CKHg2}  
  \item If $H$ is \sfc\ and $H^g \leq K$ for some $g \in G$ and some
    $p$-subgroup $K$ of $G$, then $K$ is \sfc , $Z(H^g)=C_K(H^g)$, and
    $Z(H^g) \geq Z(K)$.
    \label{lemma:CKHg3}
  \item If $Q \leq P$ and $C_P(Q)$ is a \syl p of $C_G(Q)$, then
    $QC_P(Q)$ is \sfc . \label{lemma:CKHg4}
  \end{enumerate}
  \end{lemma}
\begin{proof}
  \eqref{lemma:CKHg1} If $Z$ is a central \syl p of $C$ then $C \cong
  Z \times C/Z$ where $C/Z \cong O_{p'}C = O^pC$.

  \noindent \eqref{lemma:CKHg2} Assume that $H$ is \sfc . Then $H^g$
  is also \sfc . The $p$-subgroup $C_P(H^g)$ is contained in the
  unique \syl p $Z(H^g)$ which is contained in $H^g$. Conversely,
  assume that $H$ has property \eqref{lemma:CKHg2}. Choose $g \in
  N_G(H,P)$ so that $C_P(H^g)$ is a \syl p of $C_G(H^g)$. By
  assumption, $C_P(H^g) = C_P(H^g) \cap H^g = Z(H^g)$. This shows that
  $Z(H)^g$ is a \syl p of $C_G(H)^g$. 

  \noindent \eqref{lemma:CKHg3} Let $h \in N_G(K,P)$. Then $gh \in
  N_G(H,P)$ and $C_P(K^h) \leq C_P(H^{gh}) \leq H^{gh} \leq K^k$.
  According to \eqref{lemma:CKHg2}, $K$ is \sfc .  Since $Z(H)$ is
  central in $C_G(H)$ it is the {\em unique\/} \syl p of $C_G(H)$.
  The $p$-subgroup $C_K(H^g)$ of $C_G(H^g)$ is a subgroup of $Z(H^g)$.
  Now the chain of inclusions
  \begin{equation*}
    Z(H^g) = C_{H^g}(H^g) \leq C_K(H^g) \leq Z(H^g)
  \end{equation*}
  shows that $Z(H^g)=C_K(H^g)$. Obviously, $Z(H^g) = C_K(H^g) \geq
  C_K(K) = Z(K)$.

  \noindent \eqref{lemma:CKHg4} According to \cite[Proposition
  2.11]{puig09}, $C_P(Q)^g = C_P(Q^g)$ and $C_P((QC_P(Q))^g) \leq
  C_P(Q^g) = C_P(Q)^g \leq (QC_P(Q))^g$ for any $g \in
  N_G(QC_P(Q),G)$. Now apply item \eqref{lemma:CKHg2}.
\end{proof}


\begin{table}[t]
  \centering
  \begin{tabular}{l||l|l|l}
    $ \mathcal{C}$ &  $k^{[H]}$  &  $k_K$  &  $\chi(\mathcal{C})$  \\ \hline

 
    $\displaystyle \catc TG$   &  
    $\displaystyle \sum_{[K]} [\mu]([H],[K])$  & 
    $\displaystyle |G|^{-1}\sum_H\mu(H,K) $  &  
    $\displaystyle \sum_{[H],[K]} [\mu]([H],[K]) $  \\

    $\catc LG$  &  
    $\displaystyle \sum_{[K]} [\mu]([H],[K]) |C_G(K)|_{p'}$  &
    $\displaystyle |G|^{-1}|C_G(K)|_{p'}\sum_H\mu(H,K) $  &
    $\displaystyle \sum_{[H],[K]} [\mu]([H],[K])|C_G(K)|_{p'}$  \\
    
    $\displaystyle \catc FG$  &  
    $\displaystyle \sum_{[K]} [\mu]([H],[K]) |C_G(K)|$  &
    $\displaystyle |G|^{-1}|C_G(K)|\sum_H\mu(H,K) $  &
    $\displaystyle \sum_{[H],[K]} [\mu]([H],[K])|C_G(K)|$  \\
    
    $\displaystyle \catc OG$  &  
    $\displaystyle |H|\sum_{[K]}[\mu]([H],[K])$  &
    $\displaystyle |G|^{-1}\sum_H |H| \mu(H,K)$  &  
    $\displaystyle  \sum_{[H],[K]} |H| [\mu]([H],[K])$  \\

    $\displaystyle \cattc FG$  &  
    $\displaystyle |H| \sum_{[K]} [\mu]([H],[K]) |C_G(K)|_{p'}$  &
    $\displaystyle |G|^{-1} |C_G(K)|_{p'} \sum_H |H| \mu(H,K)$
    &
    $\displaystyle  \sum_{[H],[K]}  |H| [\mu]([H],[K])  |C_G(K)|_{p'}$ 
  \end{tabular}
  \caption{Categories of \sfc\ $p$-subgroups}
  \label{tab:catc}
\end{table}

\begin{thm}\label{thm:STLFOc}
  Weightings $k^\bullet$, co\we s $k_\bullet$, and \Euc s for the
  finite categories $\catc TG$, $\catc LG$, $\catc FG$, $\catc OG$,
  and $\cattc FG$ of $p$-selfcentralizing $p$-subgroups of $G$ are as
  in Table~\ref{tab:catc}.
\end{thm}
\begin{proof}
  This follows almost immediately from Theorem~\ref{thm:ShG}.  We
  comment on the most interesting case, $\cattc FG$.  If $H$ is \sfc\
  and $K \geq H$, then $K$ is \sfc\ and $C_K(H)$ is isomorphic to $Z(H)$
  (Lemma~\ref{lemma:CKHg}.\eqref{lemma:CKHg3}).  Equality
  \eqref{eq:tildeFG1} simplifies to
\begin{equation*}
  |C_G(H)|_{p'} |\cattc FG(H,K)| |K| = |\catn TGc(H,K)|
\end{equation*}
so that the functions
\begin{equation*}
  k^H = |G|^{-1}|H| \sum_K \mu(H,K) |C_G(K)|_{p'}, \qquad
  k_K =  |G|^{-1} |C_G(K)|_{p'} \sum_H |H| \mu(H,K)
\end{equation*}
are a \we\ and a co\we\ for $\cattc FG$ by Theorem~\ref{thm:ShG}. 
Rewriting the \we\ as
\begin{equation*}
  k^H = \frac{|H|}{|G \colon N_G(H)} \sum_{[K] \in [\catn TGc]}
  [\mu]([H],[K]) | C_G(K)|_{p'}, \qquad
  k^{[H]} = |H| \sum_{[K]} [\mu]([H],[K]) | C_G(K)|_{p'}
\end{equation*}
we calculate the \Euc\ 
\begin{equation*}
  \chi(\cattc FG) =\sum_{[H]}k^{[H]} 
  = \sum_{[H]} |H| [\mu]([H],[K]) | C_G(K)|_{p'}
\end{equation*}
as the sum of the values of the \we .
\end{proof}

The  function $[\mu]$ of Table~\ref{tab:wecat*} is
the \Mb\ function
of $[\catn TGc]$ (Proposition~\ref{prop:nupi0TG*}).  Because the \sfc\
property is upward closed (Lemma~\ref{lemma:CKHg}.\eqref{lemma:CKHg3})
the \Mb\ for $[\catn TGc]$ is simply the restriction of the \Mb\
function for $[\catn TG*]$.  Also, the {\em \we s\/} for $\catn CGc$
from Table~\ref{tab:catc} are the restrictions of the \we s for $\catn
CG*$ for $\mathcal{C}=\mathcal{T}, \mathcal{L}, \mathcal{F},
\mathcal{O}$ from Table~\ref{tab:wecat*}.

  \begin{table}[t]
    \centering 
    \begin{tabular}{c||c|c|c|c|c|c|c|c|c|c} \hline
      $n$ 
      & \multicolumn{2}{c|}{ $4-5$} & $6$ & $7$ 
      & \multicolumn{2}{c|}{$8-9$} & $10$ & $11$ 
      & \multicolumn{2}{c}{$12 - 13$} \\ \hline \hline

      $\chi(\catc O{A_n})$ 
      & \multicolumn{2}{c|}{$1/3$}
      & $1/3$ & \multicolumn{1}{c}{$2/9$} 
      & \multicolumn{2}{|c|}{$13/63$} 
      & $19/105$ & $106/945$ 
      & \multicolumn{2}{c}{$388/2835$} \\ \hline

      $\chi(\catc T{A_n})$ & \multicolumn{2}{c|}{$1/12$} &  
      $-1/24$ & $-5/72$ & 
      \multicolumn{2}{c|}{$13/4032$}  &
      $29/4480$ & $653/120960$  & 
      \multicolumn{2}{c}{$1133/1451520$}  \\ \hline

      $\chi(\catc L{A_n})$ 
      & \multicolumn{2}{c|}{$1/12$} 
      & \multicolumn{2}{c|}{$-1/24$} 
      & \multicolumn{2}{c|}{$13/4032$} 
      & \multicolumn{2}{c|}{$29/4480$}  
      & \multicolumn{2}{c}{$1133/1451520$}  \\ \hline

      $\chi(\catc F{A_n})$ 
      & \multicolumn{2}{c|}{$1/3$}  
      & \multicolumn{2}{c|}{$1/3$}  
      & \multicolumn{2}{c|}{$13/63$}  
      & \multicolumn{2}{c}{$19/105$} 
      & \multicolumn{2}{|c}{$388/2835$}       \\ \hline
      
      $\chi(\cattc F{A_n})$
      & \multicolumn{2}{c|}{$1/3$} 
      & \multicolumn{2}{c|}{$1/3$} 
      & \multicolumn{2}{c|}{$13/63$}
      & \multicolumn{2}{c|}{$19/105$}
      & \multicolumn{2}{c}{$388/2835$}
    \end{tabular}
    \caption{\Euc s at $p=2$ of \sfc\ $p$-subgroup posets and 
      categories of alternating groups computed by Magma \cite{magma}}
    \label{tab:CcAn}
  \end{table}



  We next note that the \we s for $\catn LGc$ and $\cattc FG$ can be
  computed locally (cf Proposition~\ref{prop:altwes}). 

 Fix $H$, a \sfc\ $p$-subgroup of $G$, and consider the
  projection $\catn TGc(H)= N_G(H) \to \overline{N_G(H)} = N_G(H)/H =
  \catn OGc(H)$ of the $p$-local subgroup $N_G(H)$ onto its quotient
  $N_G(H)/H$. The functor
  \begin{equation}\label{defn:OpCGfunctor}
    O^pC_G \colon \big( \catn S{\catn OGc(H)}* \big)^{\mathrm{op}} \to 
    \catn {\overline{S}}{O^pC_G(H)}{}
  \end{equation}
  takes the nonidentity $p$-subgroup $\overline{K}$ of $\catn OGc(H)$
  to the subgroup $O^pC_G(K)$ of $O^pC_G(H)$ where $K \leq N_G(H)$ is
  the preimage of $\overline{K} \leq N_G(H)/H$. For every $x \in
  O^pC_G(H)$,
  \begin{equation}\label{defn:OpCG^-1}
    O^pC_G / \langle x \rangle =
    \{ \overline{K} \in \Ob{\catn S{\catn OGc(H)}*} \mid
    O^pC_{C}(K) \ni x\}
  \end{equation}
  is the preimage of the subposet $\{ Y \mid \langle x \rangle \leq Y
  \leq O^pC_G(H) \}$ under the functor $O^pC_G$.



\begin{prop}\label{prop:wetildeFc}
  The value of the \we\ for $\cattc FG$ at the \sfc\
  $p$-subgroup $H \leq G$ is
 \begin{equation*}
    k_{\widetilde{\mathcal{F}}^c}^H =
    |G|^{-1} \sum_{x \in O^pC_G(H)} 
     \big(1- \chi(O^pC_G/\langle x \rangle) \big)
  \end{equation*}
  \end{prop}
\begin{proof}
  By Table~\ref{tab:catc} and \eqref{defn:nu} the \we\ for $\cattc FG$
  at the \sfc\ $H \leq G$ is given by
\begin{equation*}
  |G \colon H|k_{\widetilde{\mathcal{F}}^c}^H = 
  \sum_{K \in [H,N_G(H)]} \mu(H,K) |C_G(K)|_{p'} 
\end{equation*}
For any nonidentity $p$-subgroup $H \leq G$, there is a commutative
diagram
\begin{equation*}
  \xymatrix{
    Z(H) \ar[r] \ar[d] & H \ar[r] \ar[d] & H/Z(H) \ar[d] \\
    C_G(H) \ar[r] \ar[d] & {\catn TGc(H)} \ar[r] \ar[d] & {\catn
      FGc(H)} \ar[d] \\ 
    C_G(H)/Z(H) \ar[r] & {\catn OGc(H)} \ar[r] & {\cattc FG(H)} }
\end{equation*}
with exact rows and columns. 
%
%
Let $K$ be a $p$-subgroup such that $H \leq K \leq N_G(H)$. Then
$C_G(K) \leq C_G(H) \leq N_G(H)$ so that $C_G(K) = C_{N_G(H)}(K)$.  In
case $H$ is a \sfc\ subgroup of $G$, the chain of inequalities,
obtained using Lemma~\ref{lemma:CKHg}.\eqref{lemma:CKHg3},
\begin{equation*}
  Z(K) \leq Z(H) \cap C_G(K) \leq H \cap C_G(K) \leq K \cap C_G(K) =
  Z(K) 
\end{equation*}
is, in fact, a chain of identities so that $Z(K) = Z(H) \cap C_G(K) =
H \cap C_G(K)$.  The projection $\catn TGc(H)= N_G(H) \to \catn OGc(H)
= N_G(H)/H$ takes $C_G(K)$ to $\overline{C_G(K)}$ with kernel $H \cap
C_G(K) = Z(K)$, the \syl p\ of $C_G(K)$. Thus $\overline{C_G(K)} =
C_G(K)/Z(K) = O^p C_G(K)$ and $|\overline{C_G(K)}| = |C_G(K)|_{p'}$.

This means that
\begin{equation*}
  |G \colon H|k_{\widetilde{\mathcal{F}}^c}^H =
  \sum_{\overline{K} \leq \catn OGc(H)} 
  \mu(\overline{K}) |\overline{C_G(K)}| =
  |\overline{C_G(H)}| - \sum_{1 \lneqq \overline{K} \leq
    \catn OGc(H)} -\mu(\overline{K})|\overline{C_G(K)}|
\end{equation*}
We now proceed as in the proof of Proposition~\ref{prop:altwes}. The
sum in the above formula is the \Euc\ of the Grothendieck construction
on the presheaf $O^pC_G$ \eqref{defn:OpCGfunctor}. The opposite of
this Grothendieck construction is the direct sum over $x \in
O^pC_G(H)$ of the of the posets \eqref{defn:OpCG^-1}.
\end{proof}

Based on explicit computations we suspect, first, that $\chi(\catn
FGc) = \chi(\cattc FG)$ (cf Corollary~\ref{cor:chiFGeqchiFtildeG})
and, second, that the \we\ for $\cattc FG$ is supported precisely on
the \sfc\ $\catn FG{}$-radical subgroups
(Definition~\ref{defn:radical}), ie that
\begin{equation}\label{conj:FGradical}
  \sum_{x \in O^pC_G(H)} 
     \chi(O^pC_G/\langle x \rangle) 
     \neq |O^pC_G(H)| \iff O_p \cattc FG(H) = 1
\end{equation}
holds for any \sfc\ subgroup $H$ of any group $G$, cf
\eqref{conj:pradical}.  Explicit computations with Magma \cite{magma}
reveal that \eqref{conj:FGradical} is true at $p=2$ for all groups of
order $\leq 760$ and for the alternating groups $A_n$, $4 \leq n \leq
13$, of Table~\ref{tab:CcAn}.

\begin{lemma}\label{lemma:conjFGtilde}
  Let $H$ and $K$ be two nonidentity $p$-subgroups of $G$. Then $H$
  and $K$ are isomorphic in $\widetilde{\mathcal{F}}^*_G$ if and
  only if they are isomorphic in $\cat FG$.
\end{lemma}
\begin{proof}
  Suppose that $H$ and $K$ are isomorphic in
  $\widetilde{\mathcal{F}}^*_G$. Then there exist $x \in N_G(H,K)$, $y
  \in N_G(K,H)$ so that conjugation by $xy$ is an inner auto\m\ of $H$
  and conjugation by $yx$ is an inner auto\m\  of $K$. By replacing
  $y$ by another element of $yH$, if necessary, we obtain that $yx \in
  C_G(H)$. Then $yx=xy^x \in C_G(H)^x = C_G(K)$. This means that $xy$
  represents the identity of $\cat FG(H)$ and $yx$ represents the
  identity of $\cat FG(K)$.
\end{proof}

  If $P$ is a nonidentity $p$-group then
  \begin{equation*}
    \chi(\catc SP)=1, \qquad
    \chi(\catn TPc)=|P|^{-1}, \qquad
    \chi(\catn LPc)=|P|^{-1}, \qquad
     \chi(\catn FPc)= 1, \qquad
    \chi(\catn OPc)=1, \qquad
    \chi(\cattc FP)=1
  \end{equation*}
  because $\catn SPc$, $\catn OPc$, and $\cattc FP$ have $P$ as
  terminal object and $\catn LPc = \catn TPc$ is the Grothendieck
  construction for the $P$-action on $\catc SP$.
  Corollary~\ref{cor:normalsyl}.\eqref{cor:normalsyl2} shows that the
  \we\ for $\catn FP*$ and is supported on the subgroups of the form
  $C_P(x)$, $x \in P$, so that $\chi(\catn FPc) = \chi(\catn FP*) = 1$
  as these subgroups are \sfc\ by
  Lemma~\ref{lemma:CKHg}.\eqref{lemma:CKHg4}.  (By
  Example~\ref{exmp:weonnoncentric} it is not true for general groups
  $G$ that the \we\ for $\catn FG*$ is supported on the \sfc\
  subgroups and Tables~\ref{tab:FAn} and \ref{tab:CcAn} contain
  several examples of alternating groups where the nonidentity and the
  centric Frobenius categories have different \Euc s.)


\def\cprime{$'$}
\providecommand{\bysame}{\leavevmode\hbox to3em{\hrulefill}\thinspace}
\providecommand{\MR}{\relax\ifhmode\unskip\space\fi MR }
\providecommand{\MRhref}[2]{%
  \href{http://www.ams.org/mathscinet-getitem?mr=#1}{#2}
}
\providecommand{\href}[2]{#2}

\end{document}